\documentclass{article}

\usepackage{amssymb}
\usepackage{amsmath} % A mathematical package, providing elements like ``bmatrix''
\usepackage{lipsum}
\usepackage{amsfonts}
\usepackage{bbm,alltt}
\usepackage{booktabs} % Horizontal rules in tables
\usepackage{tabularx}
\usepackage{float} 
\usepackage{graphicx,epstopdf} % Required to insert images
\usepackage[caption=false]{subfig} % edmond 6/16/2017
\usepackage{blkarray}

\usepackage[colorlinks]{hyperref}
\usepackage[capitalise,noabbrev]{cleveref}	% for \cref

\usepackage{algorithm}
\usepackage{algorithmic} % Algorithm style, could alternatively use algpseudocode

\usepackage{verbatim} % edmond 6/12/2017

\newcounter{theorem}
\newtheorem{problem}[theorem]{Problem}

\ifpdf%
  \DeclareGraphicsExtensions{.eps,.pdf,.png,.jpg}
\else
  \DeclareGraphicsExtensions{.eps}
\fi

\title{An efficient method for block low-rank approximations for kernel matrix systems%
	\thanks{Version of \today.}}
\date{}
\author{
	Xin Xing%
	\thanks{School of Computational Science and Engineering, Georgia Institute of Technology,
		Atlanta, GA ({xxing33@gatech.edu}, {echow@cc.gatech.edu}).}
	\and
	Edmond Chow%
	\footnotemark[2]
}

\begin{document}

\maketitle

\begin{abstract}
In the iterative solution of dense linear systems from boundary 
integral equations or systems involving kernel matrices, the main
challenges are the expensive matrix-vector multiplication and the
storage cost which are usually tackled by hierarchical matrix
techniques such as $\mathcal{H}$ and $\mathcal{H}^2$ matrices. 
However, hierarchical matrices also have a high construction cost 
that is dominated by the low-rank approximations of the sub-blocks 
of the kernel matrix. 
%However, in hierarchical matrix construction, the low-rank
%approximations of the sub-blocks of the kernel matrix also lead
%to a high construction cost. 
In this paper, an efficient method is proposed to give a 
low-rank approximation of the kernel matrix block $K(X_0, Y_0)$ in the form of
an interpolative decomposition (ID) for a kernel function $K(x,y)$ and two properly
located point sets $X_0, Y_0$. 
The proposed method combines the ID using strong rank-revealing QR (sRRQR),
which is purely algebraic, with analytic kernel information to
reduce the construction cost of a rank-$r$ approximation from $O(r|X_0||Y_0|)$, 
for ID using sRRQR alone, to $O(r|X_0|)$ which is not related to $|Y_0|$. 
Numerical experiments show that $\mathcal{H}^2$ matrix construction
with the proposed algorithm only requires a computational cost linear in the matrix
dimension.  
\end{abstract}

%\begin{keywords}
  % 7. Keywords that describe the paper
%\end{keywords}
\section{Introduction}
%In problems such as boundary integral equations and particle simulations, the
%long-range interactions between elements or particles lead to large dense kernel
%matrices that are challenging for iterative solvers. 
In this paper, we are concerned with dense matrices generated by a translation-invariant
kernel function $K (x, y) = k(x-y)$ that satisfies the property that for any two separated 
clusters of points, $X_0 = \{ x_i \}_{i = 1}^n$ and $Y_0 = \{ y_j \}_{j = 1}^m$, 
the kernel matrix $K (X_0, Y_0) = (K (x_i, y_j))_{x_i \in X_0, y_j \in Y_0} \in \mathbb{R}^{n
\times m}$ is numerically low-rank. The low-rank property of $K
(X_0, Y_0)$ is usually evidenced by an analytic expansion with separated
variables for the kernel function, i.e.,
\begin{equation}
  K (x, y) = \sum_{i = 1}^r \psi_i (x) \phi_i (y) + R_r (x, y),  \label{expansion}
\end{equation}
where $\{\psi_i(x)\}$ and $\{\phi_i(y)\}$ are some functions of one variable. 
The remainder $R_r (x, y)$ is close to zero by requiring certain conditions on
the separation of points in $X_0$ and $Y_0$.  Such $X_0$ and $Y_0$ pairs are then said to be admissible.
Denoting the convex hulls of $X_0$ and $Y_0$ as $\mathcal{X}_0$ and $\mathcal{Y}_0$ 
respectively, the typical admissibility conditions for $X_0\times Y_0$, equivalent to 
those for $\mathcal{X}_0\times \mathcal{Y}_0$, include 
\begin{itemize}
\item strong admissibility condition:
$\min\left(\text{diam}(\mathcal{X}_0), \text{diam}(\mathcal{Y}_0)\right) \leqslant \eta\ 
\text{dist}(\mathcal{X}_0,\mathcal{Y}_0)$ for a constant $\eta$ and where 
$\text{diam}(\mathcal{X}_0)$ denotes a measure of the diameter of $\mathcal{X}_0$. 

\item weak admissibility condition: $\mathcal{X}_0\cap \mathcal{Y}_0 = \varnothing$.
\end{itemize}

For a kernel matrix with prescribed point sets, certain sub-blocks of the
kernel matrix can be associated with admissible cluster pairs and hence are
numerically low-rank. Representing these sub-blocks by various low-rank forms with
different admissibility conditions and additional constraints, hierarchical
matrix representations, such as $\mathcal{H}$ \cite{hackbusch_sparse_1999,hackbusch_sparse_2000}, 
$\mathcal{H}^2$ \cite{hackbusch_data-sparse_2002,hackbusch_$mathcalh^2$-matrices_2000}, 
HSS \cite{chandrasekaran_fast_2006} and HODLR \cite{ambikasaran_mathcalo_2013}, can
help reduce both the matrix-vector multiplication complexity and the storage cost
from $O (n^2)$ to $O (n \log^{\alpha} n)$ or even $O(n)$. Similarly, fast
matrix-vector multiplication algorithms, like FMM \cite{greengard_fast_1987,greengard_new_1997}
and panel clustering \cite{hackbusch_fast_1989}, use the same idea and are algebraically 
equivalent to certain types of the above hierarchical matrix representations.

\begin{figure}[ht]
	\centering
	\label{fig:point_set_pair}
	\subfloat[Point set pair for $\mathcal{H}$ construction]{
		\includegraphics[width=0.35\textwidth]{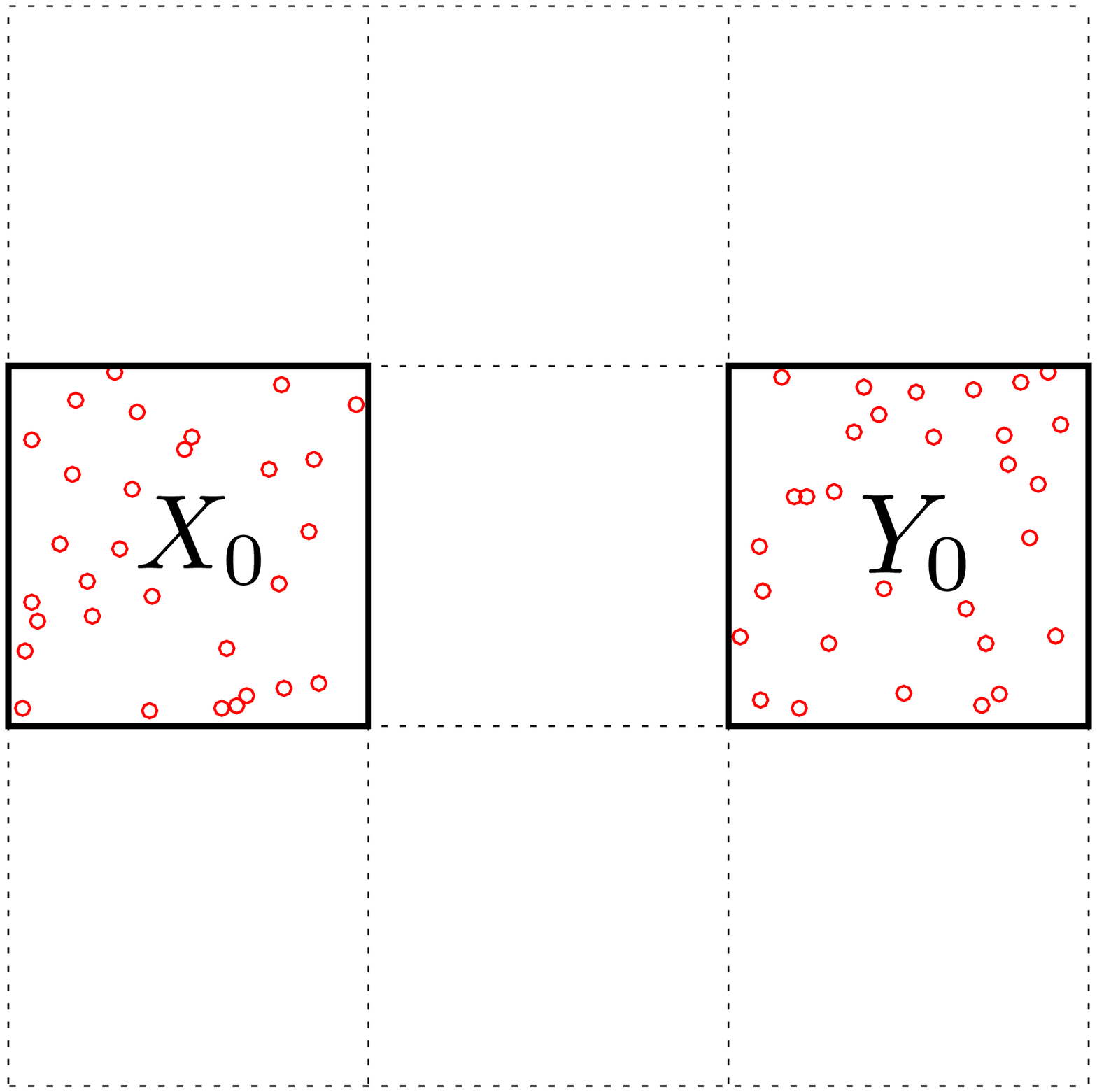}
	}
	\hspace{2em}
	\subfloat[Point set pair for $\mathcal{H}^2$ construction]{
		\includegraphics[width=0.35\textwidth]{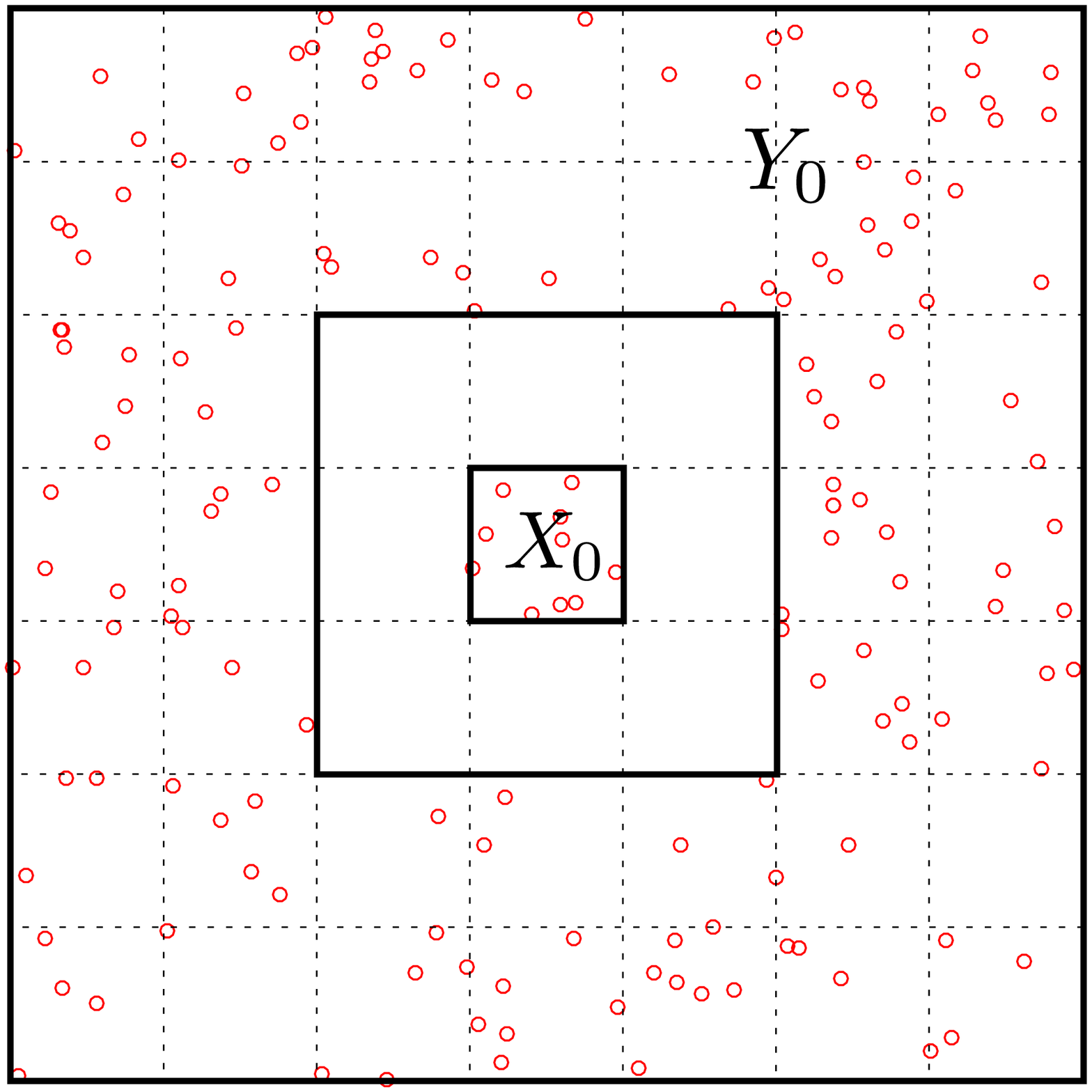}
	}
	\caption{Examples of point set pairs $X_0 \times Y_0$ in the approximated low-rank blocks with a strong 
		admissibility condition.}
\end{figure}

Although they provide great savings, hierarchical matrix representations usually have a high
construction cost that is dominated by computing the low-rank approximation of certain sub-matrices
or blocks. Specifically, the low-rank blocks approximated in $\mathcal{H}$ construction are
all of the form $K(X_0, Y_0)$ with an admissible cluster pair $X_0 \times Y_0$.
In $\mathcal{H}^2$ construction with interpolative decomposition
\cite{martinsson_fast_2011,cai_difeng_smash:_2017, ho_fast_2012} 
(referred to as ID-based $\mathcal{H}^2$ construction), the blocks are of the form
$K(X_0, Y_0)$ with $X_0$ being a cluster and $Y_0$ being the union of all clusters that are admissible with $X_0$. 
Examples in 2D of point set pairs $X_0 \times Y_0$ in both cases are shown in 
\cref{fig:point_set_pair}. It is critical to have an efficient algorithm for the low-rank 
approximation of $K(X_0, Y_0)$ with $X_0\times Y_0$ in both these cases. 

Using an algebraic approach, interpolative decomposition (ID) \cite{gu_efficient_1996,cheng_compression_2005},
QR variants, adaptive cross approximation (ACA) \cite{bebendorf_approximation_2000,bebendorf_adaptive_2003}
and randomized rank-revealing algorithms \cite{halko_finding_2011} are widely used in hierarchical
matrix construction. 
Most of these algebraic methods take at least $O (r | X_0 | | Y_0 |)$ with the obtained rank $r$. 
The only exception is ACA with complexity $O(r^2(|X_0|+|Y_0|))$ but its validity is based on 
the smoothness of the kernel function and certain admissibility conditions. 
%, and it may not be stable \cite{borm_hybrid_2005}. 

Using an analytic approach, low-rank approximations of $K(X_0, Y_0)$ can be obtained by a degenerate function
approximation of $K(x,y)$, i.e., a finite expansion with separated variables like the summation term
in \cref{expansion}. Such an approach only requires $O(r(|X_0|+|Y_0|))$ computation.
Typical strategies include Taylor expansion, as in panel clustering \cite{hackbusch_fast_1989}, and 
multipole expansion, as in FMM \cite{greengard_new_1997}.
However, the obtained rank $r$ can be much larger than those by algebraic methods and explicit 
expansions are only available for a few standard kernels.

%FMM and panel clustering essentially approximate $K (x, y)$ analytically by degenerate kernel functions, 
%i.e., a finite expansion with separated variables like the summation term in \cref{expansion}. 
%The sub-block $K (X_0, Y_0)$ can then be approximated as 
%\begin{equation}
%  K (X_0, Y_0) \approx \Psi (X_0)^T \Phi (Y_0),\quad \left[\Psi (X_0)\right]_{i j} =
%  \psi_i (x_j),\quad  \left[\Phi (Y_0)\right]_{i j}  = \phi_i(y_j), \label{eqn:analytic_approx}
%\end{equation}
%%where $\Psi (X_0)\in \mathbb{R}^{r\times n}$ and $\Phi (Y_0)\in \mathbb{R}^{r\times m}$ 
%and it only takes $O(r(|X_0|+|Y_0|))$ computation.
%%Typical strategies include Taylor expansion, interpolation and multipole expansion. 
%%The rank $r$ and certain admissibility condition are derived analytically so that the 
%%remainder function $|R_r (x, y)|$ in \cref{expansion} is bounded close to zero. 
%However, the obtained $r$ can be much larger than those by algebraic methods and explicit 
%expansions are only available for a few standard kernels. 

There are also several hybrid algebraic-analytic compression algorithms such as those used in 
kernel-independent FMM (KIFMM) \cite{ying_kernel-independent_2004}, recursive skeletonization 
\cite{martinsson_fast_2005, ho_fast_2012} and SMASH \cite{cai_difeng_smash:_2017}.
These algorithms share the same strategy of taking advantage of having an analytic
kernel but without having any explicit expansion of $K(x,y)$. This strategy is combined
with purely algebraic methods to help reduce the computational cost. However, the validity of both 
KIFMM and recursive skeletonization is only proved for kernels from potential theory, and SMASH needs
a heuristic selection of the rank for certain kernel matrix blocks and the basis functions for 
degenerate function approximation of $K(x,y)$.  

Following the same strategy as the above hybrid algorithms, we introduce a new algorithm for 
the low-rank approximation of $K(X_0,Y_0)$ for general kernel functions that implicitly
uses the putative degenerate function approximation. 
The method also uses the ID by strong rank-revealing QR (sRRQR) \cite{gu_efficient_1996} but 
reduces the construction cost from $O(r|X_0||Y_0|)$, for ID using sRRQR alone, to $O(r|X_0|)$. 
%The method also uses ID but reduces the construction cost from $O(r|X_0||Y_0|)$, for ID alone, to $O(r|X_0|)$. 
The proposed algorithm only requires kernel evaluations and can automatically determine the rank
for a given error threshold. 
\section{Background}
%For the following discussion, we only consider a symmetric kernel matrix $K(X,X)$ with a
%prescribed point set $X$ and the approximation of sub-matrices in its ID-based $\mathcal{H}^2$
%construction. The proposed algorithm can easily be adapted to $\mathcal{H}$ construction.
In this paper, we focus on the approximation of blocks in ID-based $\mathcal{H}^2$ 
construction but the same ideas can be easily adapted to $\mathcal{H}$ construction. 

Hierarchical matrix construction is based on a hierarchical partitioning of a box domain in
$\mathbb{R}^d$ where the box encloses all the prescribed points. 
Defining this box as the root level, finer partitions at subsequent levels are obtained by
recursively subdividing every box at the previous level uniformly into $2^d$ smaller boxes
until the number of points in each finest box is less than a prescribed constant. 
In ID-based $\mathcal{H}^2$ construction, 
each non-empty box $\mathcal{B}_i$ at any non-root level is associated with an ID 
approximation of a sub-block $K(X_i, Y_i)$ where $X_i$ is some subset of the points lying in 
$\mathcal{B}_i$ and $Y_i$ is some subset of the points lying in the union of boxes
at the same level that are admissible with $\mathcal{B}_i$. Readers can refer to 
\cite{ho_fast_2012, cai_difeng_smash:_2017} for more details. 

\begin{figure}[ht]
\centering
\label{fig:domain}
\subfloat[Strong admissibility condition]{
\includegraphics[width=0.33\textwidth]{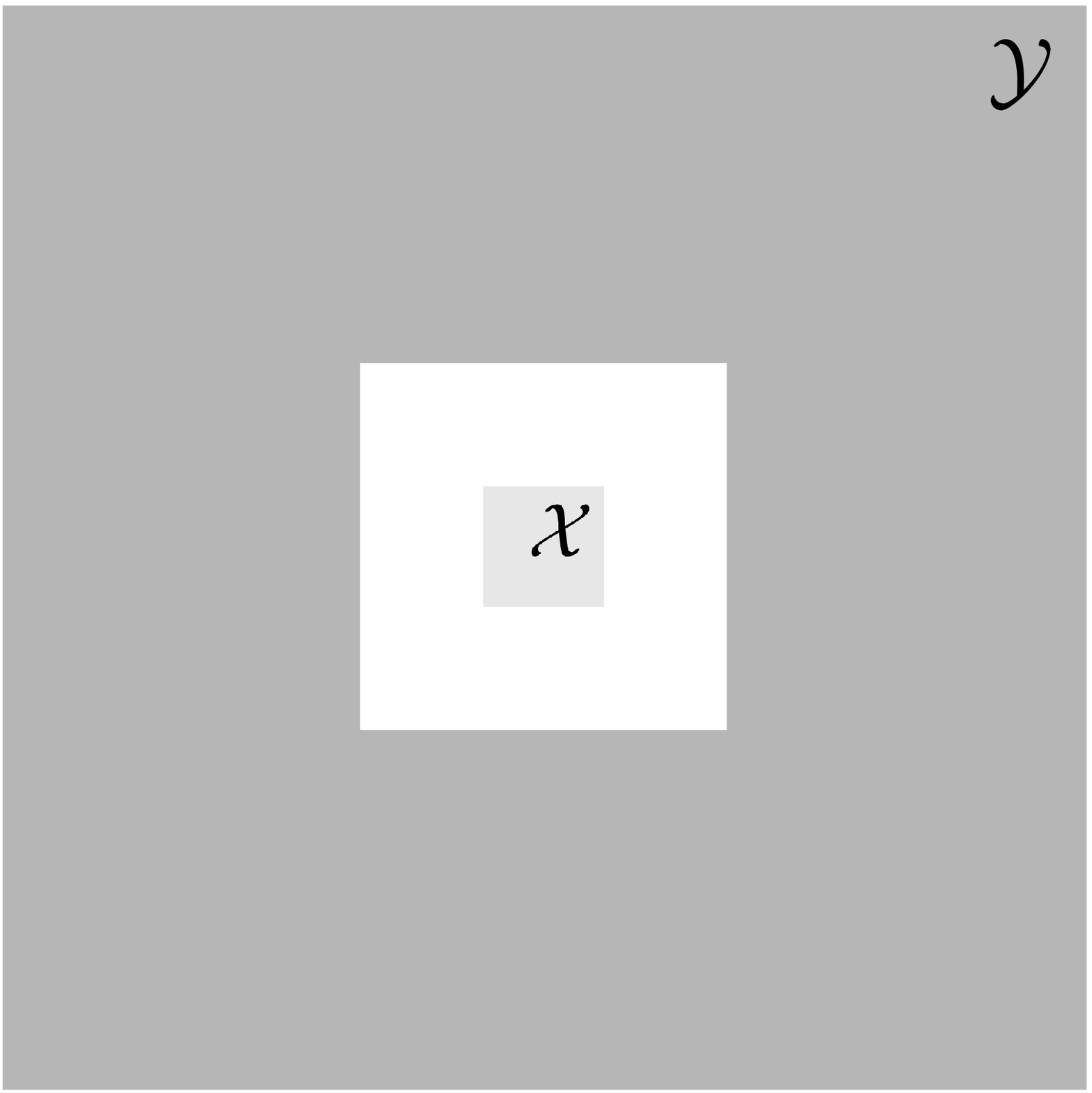}
}
\hspace{2em}
\subfloat[Weak admissibility
 condition]{
\includegraphics[width=0.33\textwidth]{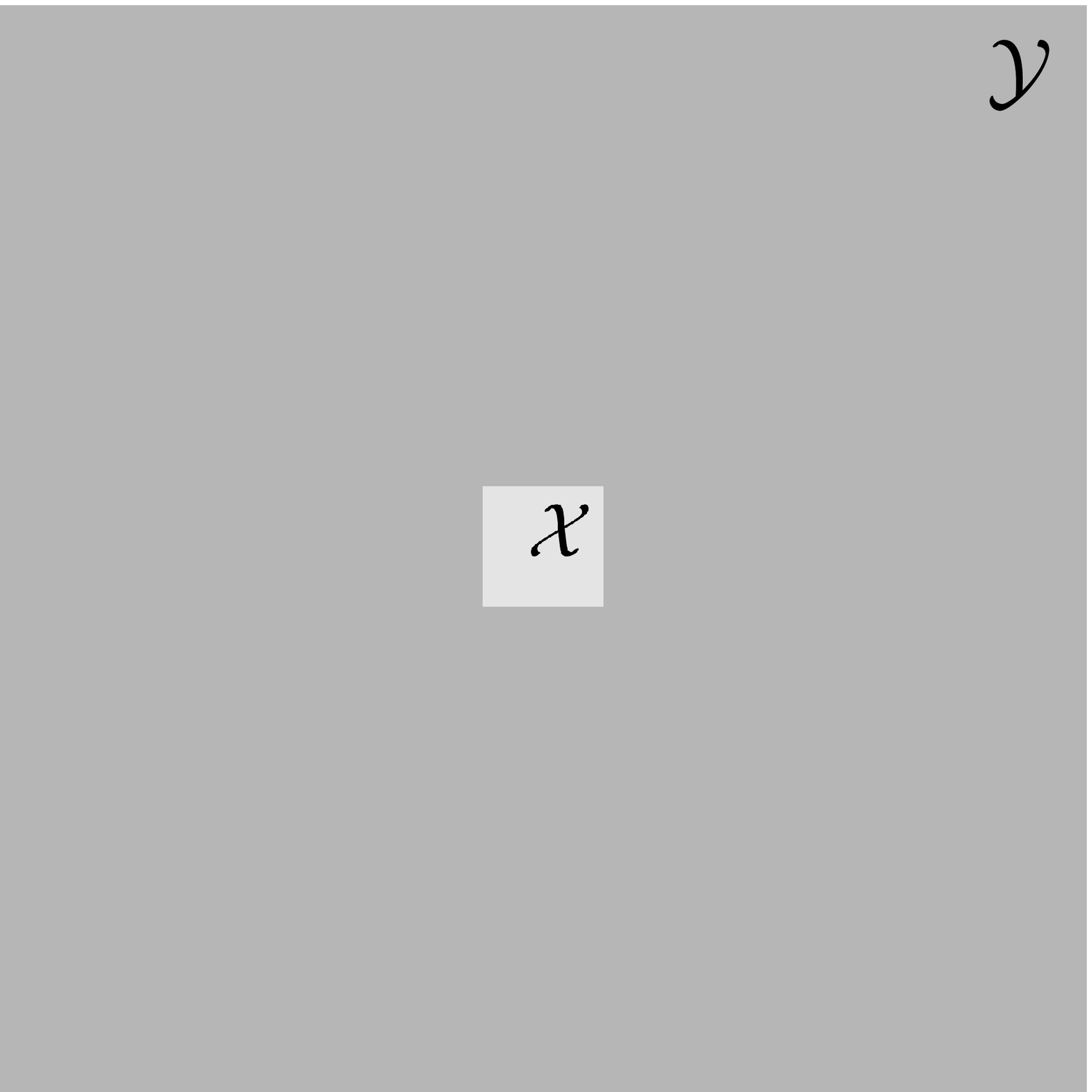}
}
\caption{Examples of the domain pair $\mathcal{X} \times \mathcal{Y}$ in 2D.}
\end{figure}

Since $K(x,y)$ is translation-invariant and boxes at the same level are of the same size, 
the approximations of $K(X_i, Y_i)$ associated with these different 
boxes $\mathcal{B}_i$ at the same level can all be unified into the following single 
problem. 
\begin{problem}\label{problem:targetmatrix}
Find an ID approximation of $K(X_0, Y_0)$ with point sets $X_0 \subset \mathcal{X}$
and $Y_0 \subset \mathcal{Y}$ where $\mathcal{X}$ is a fixed box and $\mathcal{Y}$ is
the union of all the boxes that have the same size as $\mathcal{X}$ and are admissible
with $\mathcal{X}$. 
The domain $\mathcal{Y}$ is referred to as the \textnormal{far field} of $\mathcal{X}$ by the 
prescribed admissibility condition. In practice, we only consider $\mathcal{Y}$ as a 
bounded sub-domain of the far field. 
Examples of $\mathcal{X}\times \mathcal{Y}$ are illustrated in \cref{fig:domain}.
%Find an ID approximation of $K(X_0, Y_0)$ with point sets $X_0 \subset \mathcal{X}$
%and $Y_0 \subset \mathcal{Y}$ where $\mathcal{X}$ is a fixed box and $\mathcal{Y}$ is
%a bounded sub-domain of the ``far field'' of $\mathcal{X}$. The \textnormal{far field} of 
%$\mathcal{X}$ by an admissibility condition is defined as the union of all the boxes
%that have the same size as $\mathcal{X}$ and are admissible with $\mathcal{X}$. 
%Examples of $\mathcal{X}\times \mathcal{Y}$ in 2D are illustrated in \cref{fig:domain}.
\end{problem}

\begin{figure}[ht]
	\centering
	\label{fig:illustration_conversion}
	\includegraphics[width=0.9\textwidth]{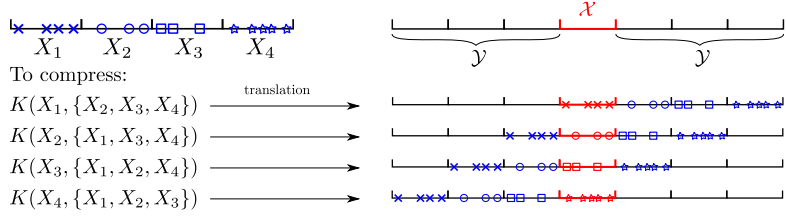}
	\caption{Example of low-rank approximations in one partition level of the ID-based 
	$\mathcal{H}^2$ construction. The weak admissibility condition is applied but the strong admissibility scenario 
	can be handled similarly. $X_1$, $X_2$, $X_3$ and $X_4$ are point sets in the subintervals 
	at the finer level of a two-level hierarchical partitioning of the whole interval. 
	For each point set $X_i$, the ID approximation of $K(X_i, \cup_{j=1}^4X_j\backslash X_i)$
	is needed. With a proper selection of domain pair $\mathcal{X}\times\mathcal{Y}$ as shown in the 
	figure, every approximation problem can be converted to \cref{problem:targetmatrix} using translation.
	}
\end{figure}

A simple one-dimensional example in \cref{fig:illustration_conversion} illustrates the associated low-rank
approximations needed in one level of ID-based $\mathcal{H}^2$ construction and
the way to convert these approximation problems into \cref{problem:targetmatrix} using translations. 

The interpolative decomposition \cite{gu_efficient_1996,cheng_compression_2005} is 
extensively used in this paper. Since somewhat different definitions exist in the 
literature, we give our definition as follows.
Given a matrix $A\in \mathbb{R}^{n\times m}$, a rank-$k$ ID approximation of $A$ is 
$U A_J$ where $A_J \in \mathbb{R}^{k\times m}$ is a row subset of $A$ and entries of 
$U \in \mathbb{R}^{n\times k}$ are bounded. ID here only refers to a form of 
decomposition and there exist many ways of computing it with different accuracies. 
In particular, minimizing the approximation error in the Frobenius norm, 
the optimal $U$ for an ID with a fixed $A_J$ can be calculated as $U = A A_J^{\dag}$ 
by projecting each row of $A$ onto the row space of $A_J$. 
%%Following \cite{ho_fast_2012}, we call $U$ and $A_J$ the projection and repleton
%%matrices, respectively. 

Define $UA_J$ as an ID with error threshold $\varepsilon$ if the norm of each row of the
error matrix $A - U A_J$ is bounded by $\varepsilon$.
The ID can be calculated algebraically by applying strong rank-revealing QR (sRRQR) 
\cite{gu_efficient_1996} to $A^T$. Truncating the obtained sRRQR decomposition of $A^T$
with absolute error threshold $\varepsilon$ as
\[
A^T P = \left(A_1^T\ A_2^T\right) 
	  = \left(Q_1 \ Q_2\right) \left(\begin{smallmatrix}
	  R_{11} & R_{12} \\  & R_{22}
	  \end{smallmatrix}\right)
\approx Q_1 \left(R_{11} \ R_{12}\right)
		=  A_1^T \left(I \ \ R_{11}^{-1}R_{12}\right),
\]
the ID with error threshold $\varepsilon$ is then 
$A \approx P  \left(\begin{smallmatrix}
I \\
(R_{11}^{-1}R_{12})^T
\end{smallmatrix}\right) A_1$ where $P$ is a permutation matrix. sRRQR can guarantee that
the entries of $R_{11}^{-1}R_{12}$ are bounded by a pre-specified parameter $C \geqslant 1$. 
The complexity of this algorithm is typically $O(rnm)$ but, in rare cases, it
may become $O(n^2m)$. 
%For the following discussion, an ID approximation always refers to 
%the one calculated by sRRQR unless otherwise specified. 

\subsection{Accelerated compression via a proxy surface\label{sec:surface}} 
\begin{figure}[ht]
\centering
\label{fig:proxy_surface}
\includegraphics[width=0.75\textwidth]{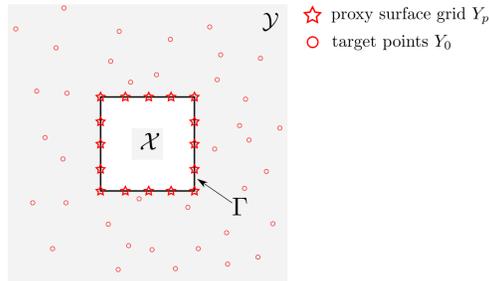}
\caption{Accelerated compression via a proxy surface. The matrix to be directly compressed changes 
from the target matrix $K(X_0, Y_0)$ to matrix $K(X_0, Y_p)$ with a constant column size $|Y_p|$, 
regardless of how many points $Y_0$ there are in $\mathcal{Y}$.}
\end{figure}

For \cref{problem:targetmatrix} with kernels from potential theory, Martinsson 
and Rokhlin \cite{martinsson_fast_2005} accelerate the ID approximation by using the concept 
of a proxy surface. Specifically, take the Laplace kernel $K(x,y)$ in 2D as an example and
consider the domain pair $\mathcal{X} \times \mathcal{Y}$ and the interior boundary of 
$\mathcal{Y}$, denoted as $\Gamma$, shown in \cref{fig:proxy_surface}. 

By virtue of Green's Theorem, the potential field in 
$\mathcal{X}$ generated by charges at $Y_0 \subset \mathcal{Y}$ can be equivalently generated by charges on 
$\Gamma$ which encloses $\mathcal{X}$. The surface $\Gamma$ is referred to as a $\textit{proxy surface}$ 
in \cite{ho_fast_2012}. 
Discretizing $\Gamma$ with a grid point set $Y_p$, it can be proved \cite{martinsson_fast_2005} that 
\begin{equation}\label{eqn:proxy_apprx}
K(x, Y_0) \approx K(x, Y_p) W, \quad \forall x \in \mathcal{X},
\end{equation}
where $W$ is a discrete approximation of the operator that maps charges at $Y_0$ 
to an equivalent charge distribution on $\Gamma$ and $\|W\|_2$ is 
bounded as a consequence of Green's Theorem. Thus, the target matrix $K(X_0, Y_0)$ can
be approximated as $K(X_0, Y_p) W$ and compressing $K(X_0, Y_0)$ directly by ID using sRRQR
is accelerated as follows. 

First find an ID approximation of $K(X_0, Y_p)$ as $UK(X_\text{rep},Y_p)$ by sRRQR where
$X_\text{rep}\subset X_0$ denotes the ``representative'' point subset associated with the selected row subset
and $U$ is the matrix obtained from the ID.  
The approximation of $K(X_0, Y_0)$ is then defined as 
$
U K(X_\text{rep}, Y_0)
$. By \cref{eqn:proxy_apprx}, the approximation error can be bounded as
\begin{align*}
\|K(X_0, Y_0) - U K(X_\text{rep}, Y_0)\|_F 
& \approx \|\left(K(X_0, Y_p) - U K(X_\text{rep}, Y_p)\right)W\|_F \\
& \leqslant \|K(X_0, Y_p) - U K(X_\text{rep}, Y_p)\|_F \|W\|_2,
\end{align*}
and thus the accelerated approximation can control the error. 

The number of points to discretize $\Gamma$ (i.e., $|Y_p|$) is heuristically decided and
only depends on the desired precision and the geometry of $\mathcal{X}$ and $\Gamma$. 
Thus, the algorithm complexity, i.e., $O(|X_\text{rep}||X_0||Y_p|)$, is independent of $|Y_0|$.  
Practically, this method only applies to $\mathcal{X} \times \mathcal{Y}$ with strong admissibility 
conditions since when $\mathcal{X}$ and $\Gamma$ are close, very large $|Y_p|$ is needed due 
to the singularity of $K(x,y)$.

The key for this method is the relation \cref{eqn:proxy_apprx} and the well-conditioning
of $W$ that are both analytically derived from Green's Theorem. Thus, the method is only rigorously
valid for kernels from potential theory and its generalization to certain problems may deteriorate
or require further modifications as discussed in \cite{martinsson_fast_2005}. 
The above method will be referred to as {\it proxy-surface method}.
 
In this paper, we develop an analogous compression algorithm for general kernels that only depends
on the existence of an accurate degenerate function 
approximation of $K(x,y)$.
In the new algorithm, the above heuristically selected point set $Y_p$ that discretizes $\Gamma$ 
will instead be selected from the whole domain $\mathcal{Y}$.  

\subsection{Notation}	
Let $\mathcal{X} \times \mathcal{Y}$ be a compact domain pair in $\mathbb{R}^d$ as described in 
\cref{problem:targetmatrix} with a certain admissibility condition and $K(x,y)$ be a 
translation-invariant kernel function in $\mathbb{R}^d$  and smooth in $\mathcal{X} \times \mathcal{Y}$. 
Denote the Karhunen-Lo{\`e}ve (KL) decomposition of $K(x,y)$ over $\mathcal{X} \times \mathcal{Y}$ as
\begin{equation}
  K (x, y) = \sum_{i = 1}^{\infty} \sigma_i \psi_i (x) \phi_i (y), 
  \quad x \in \mathcal{X}, y \in \mathcal{Y},
\end{equation}  
where $\{\psi_i(x)\}$ and $\{\phi_i(y)\}$ are sets of orthonormal functions in $\mathcal{X}$ and
$\mathcal{Y}$ respectively and $\{\sigma_i\}$ is a sequence of decaying non-negative real numbers. 
As the series converges uniformly, there exists a minimal index $r$ such that
\begin{equation}
  \left| K (x, y) - \sum_{i = 1}^r \sigma_i \psi_i (x) \phi_i (y) \right|
  \leqslant \varepsilon_{\text{machine}},\quad \forall x \in \mathcal{X},\ \forall y \in \mathcal{Y}.
  \label{ineqn:truncation}
\end{equation}

For any finite point sets $X_0 \subset \mathcal{X}$ and $Y_0 \subset \mathcal{Y}$, 
$K (X_0, Y_0)$ can then be written as
\begin{equation}\label{eqn:exact_rep}
  K (X_0, Y_0) = \Psi (X_0)^T 
  \left(\begin{smallmatrix}
  \sigma_1 &  &  & \\
  & \sigma_2 &  & \\
  & & \ldots &  & \\
  & & & \sigma_r
  \end{smallmatrix}\right)  
  \Phi (Y_0) + E, 
\end{equation}
where the entries of $E$ are bounded by $\varepsilon_{\text{machine}}$, $\Psi(x)=(\psi_i (x))_{i = 1}^r$
and $\Phi(y)= (\phi_i (y))_{i = 1}^r$ are column vector functions and $\Psi(X_0)$
with $X_0 = \{x_i\}_{i=1}^n$ is defined as
\begin{equation}\label{eqn:def_func_set}
  \Psi ( X_0 ) = 
  \left( \Psi (x_1) \ \Psi (x_2) \ \ldots \ \Psi (x_n) \right) 
  \in \mathbb{R}^{r \times n}. 
\end{equation}
We define $\Phi(Y_0)$ in the same manner. 
In this paper, the evaluation of any function over a point set is defined in the same way as
above. In particular, for a scalar function like $\psi_i(x)$, $\psi_i(X_0)$ denotes a row vector
of length $|X_0|$. Based on \cref{eqn:exact_rep}, the numerical rank of $K (X_0, Y_0)$ for any
point set pair $X_0\times Y_0$ is $r$ or less. 

For the following discussion, we consider the simplified case where $K(x,y)$ over $\mathcal{X} \times \mathcal{Y}$
is a degenerate function, i.e., its KL expansion only has a finite number of terms as
\begin{equation}\label{eqn:truncation}
  K (x, y) = \sum_{i = 1}^{r_\text{KL}} \sigma_i \psi_i (x) \phi_i (y), 
  \quad x \in \mathcal{X}, y \in \mathcal{Y}, r_\text{KL} < \infty.
\end{equation} 
A non-degenerate kernel over $\mathcal{X} \times \mathcal{Y}$ can be approximated by the first $r$ 
terms of its KL expansion with $r$ satisfying \cref{ineqn:truncation}. 
The effect of the error $\sum_{i=r+1}^{\infty} \sigma_i \psi_i(x) \phi_i(y)$, which is bounded by $\varepsilon_\text{machine}$, 
on the following proposed algorithm can be analyzed through a stability analysis which we leave as future work.

\section{Algorithm description}
Denoting $X_0\subset \mathcal{X}$ and $Y_0\subset \mathcal{Y}$ as given point
sets, our goal is to find an approximation of the target matrix $K(X_0, Y_0)$
in the form of an ID that is more efficient than using sRRQR. 

The key for ID approximation is to find a row subset of $K(X_0,Y_0)$, i.e., a subset of $\{K(x_i, Y_0)\}_{x_i \in X_0}$, 
whose span is close to each row vector $K(x_i, Y_0)$. 
Regarding $K(x_i,Y_0)$ as function values of $K(x_i,y)$ at $Y_0$, it is then sufficient to consider the above problem 
in terms of the functions $\{K(x_i, y)\}_{x_i\in X_0}$ in the domain $\mathcal{Y}$. 
Specifically, we seek a function subset of $\{K(x_i, y)\}_{x_i\in X_0}$ whose span is close to each function 
$K(x_i, y)$ in $\mathcal{Y}$.

The above ``ID approximation'' of functions $\{K(x_i, y)\}_{x_i \in X_0}$ in $\mathcal{Y}$ is a continuous problem. 
Heuristically, we can use a uniform grid point set $Y_p$ in $\mathcal{Y}$ to discretize the function 
$K(x_i,y)$ as $K(x_i, Y_p)$, transforming the problem into finding an ID approximation of $K(X_0, Y_p)$. 
In general, $Y_p$ should be dense enough to accurately characterize $K(x_i,y)$ in $\mathcal{Y}$ 
but this is inefficient in general. 
However, by the finite KL expansion in \cref{eqn:truncation}, $K(x_i, y)$ for any $x_i \in X_0$ is 
in the $r_\text{KL}$-dimensional function space spanned by $\phi_1(y)$, $\phi_2(y)$, \ldots, $\phi_{r_\text{KL}}(y)$.
Since $\{\phi_i(y)\}_{i=1}^{r_\text{KL}}$ are orthonormal in $\mathcal{Y}$, to 
uniquely determine $K(x_i, y)$, the selected finite point set $Y_p \subset \mathcal{Y}$
only needs to satisfy
\begin{equation} \label{prop:Ysample}
| Y_p | \geqslant r_\text{KL},\  \text{col} (\Phi (Y_p)) =\mathbb{R}^{r_\text{KL}}.
\end{equation}
Importantly, these are points selected from the entire domain $\mathcal{Y}$, 
not just from $Y_0$ or the boundary $\Gamma$ of $\mathcal{Y}$ and we refer 
these points as \textit{proxy points}. 
We expect an effective $Y_p$ to satisfy $|Y_p|\sim O(r_\text{KL})$ and $|Y_p| \ll |Y_0|$ in real situations.

With a proxy point set $Y_p$ that satisfies \cref{prop:Ysample}, the following algorithm 
is proposed to find an ID approximation of $K(X_0, Y_0)$ through an ``ID approximation'' of 
$\{K(x_i, y)\}_{x_i \in X_0}$ in $\mathcal{Y}$. 

\paragraph{Step 1} Find an ID approximation of $K(X_0, Y_p)$ with error threshold $\varepsilon$ by sRRQR as 
\begin{equation} \label{eqn:sampleID}
K(X_0, Y_p) \approx W_\text{rep} K(X_\text{rep}, Y_p),
\end{equation}
where $X_\text{rep} \subset X_0$ denotes the point set associated with the selected
row subset and $W_\text{rep}=K(X_0,Y_p)K(X_\text{rep},Y_p)^\dag$ is the obtained matrix from the ID. 

\paragraph{Step 2} For each $x_i \in X_0$, denote the $i$th row of $W_\text{rep}$ 
as $w_i$ and approximate the function $K (x_i, y)$ as 
\begin{equation}\label{apprx:function}
K(x_i, y) \approx w_i^T K (X_{\text{rep}}, y), \quad y \in \mathcal{Y}.
\end{equation}
It is expected that each $K(x_i, y)$ is close to the span of $\{K(x_j, y)\}_{x_j\in X_\text{rep}}$ 
and the associated approximation above has small error. 
Evaluating the functions in \cref{apprx:function} at $Y_0$, row vector $K(x_i, Y_0)$ can be approximated
by $w_i^TK(X_\text{rep}, Y_0)$ and a rank-$|X_\text{rep}|$ ID approximation is then defined as
\begin{equation} \label{eqn:lowrankapprx}
  K (X_0, Y_0) \approx W_{\text{rep}} K(X_{\text{rep}}, Y_0).
\end{equation}
Both $W_\text{rep}$ and $X_{\text{rep}}$ are calculated in \textit{Step 1} and only require 
$O(|X_\text{rep}||X_0||Y_p|)$ computation which is independent of $|Y_0|$. 

To summarize, the proposed algorithm calculates $\{w_i\}$ and $X_\text{rep}$ to minimize the
function approximation error $K(x_i, y)-w_i^TK(X_\text{rep}, y)$ at $Y_p$ to help make the error
small over the whole domain $\mathcal{Y}$. Thus, for any $Y_0\subset \mathcal{Y}$,
the proposed approximation \cref{eqn:lowrankapprx} has its entry-wise error bounded by 
$\max_{x_i\in X_0} \|K(x_i,y)-w_i^TK(X_\text{rep},y)\|_\infty$. 
It is worth noting that the rank $|X_\text{rep}|$ is fixed for any $Y_0\subset\mathcal{Y}$ and 
is only related to $X_0$ and the error threshold $\varepsilon$.   
A better ID approximation with the selected $X_\text{rep}$ can be obtained by replacing $W_{\text{rep}}$ with
 $K(X_0,Y_0) K(X_{\text{rep}}, Y_0)^{\dag}$ but this requires a computational cost linear in $|Y_0|$.
Also, $K(X_\text{rep}, Y_0)$ does not need to be explicitly formulated in ID-based 
$\mathcal{H}^2$ construction which avoids $|Y_0|$-dependent calculations.

\section{Algorithm analysis}
For any $x \in \mathcal{X}$, define $g_{x} (y) = K (x, y)$ as a function of $y$ in  $\mathcal{Y}$. 
%similarly as in \cref{eqn:def_func_set}. 
By the finite KL expansion, $g_x(y)$ can be represented as
\begin{equation} \label{eqn:gfunc}
    g_x (y) = u_x^T \Phi (y), \quad u_x = (\sigma_i \psi_i (x))_{i=1}^{r_\text{KL}}, \quad \Phi(y) = (\phi_i(y))_{i=1}^{r_\text{KL}}.
\end{equation}
With a proxy point set $Y_p$ that satisfies \cref{prop:Ysample}, substitute $Y_p$ for  $y$ in \cref{eqn:gfunc} 
and solve for $u_{x}$ as $u_x^T = g_x (Y_p) \Phi (Y_p)^{\dag}$ where $g_x(Y_p)$ is a row vector defined in the obvious way.
Thus, $g_x(y)$ can be represented in terms of $g_x(Y_p)$ as
\begin{equation}\label{eqn:representation}
  g_x (y) = g_x (Y_p) \Phi (Y_p)^{\dag} \Phi (y),\quad  \forall y \in \mathcal{Y}.
\end{equation}

We can then estimate the error of both the function approximation \cref{apprx:function} and the 
proposed ID approximation \cref{eqn:lowrankapprx}. 
Denote the error of \cref{apprx:function} for each $x_i \in X_0$ as
\begin{equation}\label{eqn:errdef}
  e_i (y) = K (x_i, y) - w_i^T K (X_{\text{rep}}, y),\quad y \in \mathcal{Y}.
\end{equation}
The error vector of the $i$th row approximation in \cref{eqn:sampleID} and \cref{eqn:lowrankapprx} 
can be exactly represented as $e_i(Y_p)$ and $e_i(Y_0)$, respectively. 
Since the error threshold for the ID approximation of $K(X_0, Y_p)$ is $\varepsilon$, 
each error vector $e_i (Y_p)$ satisfies $\| e_i (Y_p) \|_2 < \varepsilon$.

Note that $e_i (y)$ is a linear combination of $g_{x_i}(y)$ and $\{g_{x_j}(y)\}_{x_j \in X_\text{rep}}$ 
and \cref{eqn:representation} holds for $g_x (y)$ with any $x\in \mathcal{X}$. 
Thus, $e_i(y)$ has the similar representation
\begin{equation}
  e_i (y) = e_i (Y_p) \Phi (Y_p)^{\dag} \Phi (y), \label{eqn:error_representation}
\end{equation}
and $e_i(y)$ and $e_i(Y_0)$ can be bounded by the error at $Y_p$ as
\begin{align}
  | e_i (y) | &\leqslant \| e_i (Y_p) \|_2 \| \Phi (Y_p)^{\dag} \Phi (y) \|_2
  			  \leqslant \varepsilon \| \Phi (Y_p)^{\dag} \Phi (y) \|_2,
  			   \label{ineqn:errorbound1} \\
  \| e_i (Y_0) \|_2 &\leqslant \| e_i (Y_p) \|_2 \| \Phi (Y_p)^{\dag} \Phi (Y_0) \|_2
			  \leqslant \varepsilon  \| \Phi(Y_p)^{\dag} \Phi(Y_0) \|_2.
  			   \label{ineqn:errorbound2}		  
\end{align}
If the choice of $Y_p$ can guarantee $\|\Phi(Y_p)^{\dag}\Phi(Y_0)\|_2$ to be small, 
the proposed approximation can then be good and its error can be controlled by $\varepsilon$. 
\section{Selection of the proxy point set}\label{sec:proxy_point}

In the proposed algorithm above, the choice of $Y_p$ is flexible but critical. 
The only requirement for $Y_p$ is the condition \cref{prop:Ysample} and we desire that 
$\|\Phi(Y_p)^\dag\Phi(Y_0)\|_2$ is small. 

By the continuity of functions $\{\phi_i(y)\}_{i=1}^{r_\text{KL}}$ in the compact domain $\mathcal{Y}$, 
$\Phi(y)$ is bounded and thus $\|\Phi(Y_p)^\dag\Phi(Y_0)\|_2$ is also bounded for any $Y_p$ satisfying
\cref{prop:Ysample}. 
The number of points in $Y_p$ also matters since adding more points to $Y_p$ can reduce 
$\|\Phi(Y_p)^\dagger \Phi(Y_0)\|_2$ monotonically but larger $|Y_p|$ leads to more computation for 
the ID approximation of $K(X_0, Y_p)$. Thus, a constraint like $| Y_p | = O (r_\text{KL})$
is necessary to balance the trade-off.  

However, $\psi_i(x)$, $\phi_i(y)$ and $r_\text{KL}$ are usually not available for a general kernel 
function over a domain pair.
We only assume that $K(x,y)$ over $\mathcal{X}\times \mathcal{Y}$ has a finite KL expansion. 
For a non-degenerate kernel, this assumption means that there is a truncation of the KL expansion 
with error satisfying \cref{ineqn:truncation}. In both degenerate and non-degenerate cases, the number 
of expansion terms (i.e., $r_\text{KL}$) only depends on the kernel and the domain pair $\mathcal{X}\times\mathcal{Y}$. 
Thus, condition \cref{prop:Ysample} cannot be directly checked for any point set $Y_p$. 
The only property we can use is based on \cref{eqn:exact_rep}, that  
$r_\text{KL}$ is an upper bound of the numerical rank of $K(X_0, Y_0)$ 
with any point set pair $X_0\times Y_0 \subset \mathcal{X} \times \mathcal{Y}$. 

Based on the analysis above, the first method to select $Y_p$ is proposed as follows.

\paragraph{Random Selection} 
Choose $Y_p$ as a set of points that are randomly and uniformly distributed in $\mathcal{Y}$ so that condition
\cref{prop:Ysample} is likely to hold. The size of $Y_p$ can be heuristically decided as the maximum 
numerical rank of $K(X_0, Y_0)$, with some tentative $X_0\times Y_0 \subset \mathcal{X} \times \mathcal{Y}$, 
plus a small redundancy constant. 

This selection turns out to be effective in many numerical tests. However, it does not 
guarantee the scaling factor $\|\Phi(Y_p)^\dag\Phi(Y_0)\|_2$ to be small and thus the proposed 
algorithm may have much larger error than that of algebraic methods with the same approximation 
rank in some cases. 

A better selection of $Y_p$ should also try to minimize $\| \Phi (Y_p)^{\dag} \Phi (Y_0) \|_2$. 
Since the point distribution $Y_0$ is problem-dependent and $\| \Phi (Y_p)^{\dag} \Phi (Y_0) \|_2 $
can be bounded as
\[
  \| \Phi (Y_p)^{\dag} \Phi (Y_0) \|_2 
  \leqslant \| \Phi (Y_p)^{\dag} \Phi (Y_0) \|_F 
  \leqslant \sqrt{| Y_0 |} \max_{y \in \mathcal{Y}} \| \Phi (Y_p)^{\dag} \Phi (y) \|_2,
\]
we only need to consider the choice of $Y_p$ to minimize 
$\| \Phi (Y_p)^{\dag} \Phi (y) \|_2$. 

Denote $S_{Y_p} (y) = \Phi (Y_p)^{\dag} \Phi (y)$ as the solution of the least squares
problem 
\begin{equation} \label{eqn:LS1}
\Phi (Y_p) S_{Y_p} (y) = \Phi (y). 
\end{equation}
Since $\{\psi_i(x)\}_{i=1}^{r_{KL}}$ are orthonormal in $\mathcal{X}$ and $u_x = (\sigma_i \psi_i(x))_{i=1}^{r_\text{KL}}$, 
we can select a point set $X_p \subset \mathcal{X}$ such that 
\begin{equation}\label{def:Xsample}
| X_p | \geqslant {r_\text{KL}},\ \text{span} \{ u_x, x\in X_p\} = \mathbb{R}^{r_\text{KL}} = \text{col}(\Psi(X_p)) . 
\end{equation}
Based on $K(x,y) = u_x^T \Phi(y)$ and \cref{eqn:LS1}, it can be verified that $S_{Y_p}(y)$
is also the solution of the least squares problem $K (X_p, Y_p) S_{Y_p} (y) = K (X_p, y)$.
As a result, $S_{Y_p}(y)$ can be represented differently as
\begin{equation}\label{eqn:SYp}
S_{Y_p}(y) = K(X_p, Y_p)^{\dag} K(X_p, y) = \Phi(Y_p)\Phi(y),
\end{equation}
with any $X_p$ satisfying \cref{def:Xsample}. By this new representation, the second method for the 
selection of $Y_p$ is described as follows.

\paragraph{ID Selection}
Select two point sets $X_d \subset \mathcal{X}$ and $Y_d \subset \mathcal{Y}$ that are
dense enough so that \cref{def:Xsample} and \cref{prop:Ysample} are likely to hold. 
Find an ID approximation of $K(X_d, Y_d)^T$ by sRRQR as
\[
K(X_d, Y_d) \approx K(X_d, Y_\text{rep}) U 
	= K(X_d, Y_\text{rep}) \Big(K(X_d, Y_\text{rep})^\dag K(X_d, Y_d)\Big)
\]
where the error threshold is set as $\sqrt{|X_d|}\varepsilon_{\text{machine}}$ to 
estimate the numerical rank of $K(X_d, Y_d)$ and it is expected that 
the rank $|Y_\text{rep}|\approx r_\text{KL}$. Then, the point subset $Y_\text{rep}$ selected by the 
ID approximation can be used as the proxy point set $Y_p$. The detailed algorithm is shown in \cref{alg:proxy}.

\begin{algorithm}
\caption{Selection of the proxy point set $Y_p$}
\label{alg:proxy}
\begin{algorithmic}[1]
\REQUIRE kernel function $K(x,y)$, domain pair $\mathcal{X} \times \mathcal{Y}$.
\ENSURE proxy point set $Y_p$. 

\STATE Select uniform grid point sets $X_d$ and $Y_d$ in $\mathcal{X}$ 
and $\mathcal{Y}$ respectively. $Y_d$ is referred to as a candidate set.
\STATE Find an ID approximation of $K(X_d, Y_d)^T$ with error threshold
 $\sqrt{|X_d|}\varepsilon_\text{machine}$ by sRRQR.
 Define $Y_p$ as the point subset of $Y_d$ associated with 
 the selected rows.
\STATE If $|Y_p| = \min(|X_d|, |Y_d|)$, select denser grids for $X_d$, $Y_d$ and repeat
 the whole process. Otherwise, $Y_p$ is the selected proxy point set.
\end{algorithmic}
\end{algorithm}

By the sRRQR used in the ID approximation, the entries of $U$ are bounded by a pre-specified parameter $C\geqslant 1$. 
Equivalently, for any $y \in Y_d$, $K(X_d, Y_\text{rep})^\dag K(X_d, y)$, as a column of $U$, has all its 
entries bounded by $C$. Since $K(x,y)$ is smooth and $Y_d$ is dense in $\mathcal{Y}$, entries 
of $K(X_d, Y_\text{rep})^{\dag}K(X_d, y)$ can also be roughly bounded by $C$ for any $y \in \mathcal{Y}$.  
Thus, it holds that $\|S_{Y_p}(y)\|_2 \lesssim \sqrt{| Y_p |}C$ for any 
$y \in \mathcal{Y}$. By the inequalities \cref{ineqn:errorbound1} and \cref{ineqn:errorbound2}, we can 
obtain upper bounds for $|e_i(y)|$ and $\|e_i(Y_0)\|_2$ as
\begin{align}
	|e_i(y)|  & \lesssim \sqrt{|Y_p|} C \varepsilon, \qquad \forall y \in \mathcal{Y}, \label{ineqn:bound}\\
	\|e_i(Y_0)\|_2 & \lesssim \sqrt{|Y_0||Y_p|} C \varepsilon. \label{ineqn:rowerrorbound}
\end{align}
Since both $|Y_p| \sim O(r_\text{KL})$ and $C$ are small, we expect the average entry-wise 
error $\|e_i(Y_0)\|_2/\sqrt{|Y_0|}$ to be $O(\varepsilon)$.

\cref{alg:proxy} has complexity $O(|Y_p||X_d||Y_d|)$ where $|Y_p|$, as an estimate of $r_\text{KL}$, only depends 
on the kernel and the domain pair $\mathcal{X}\times \mathcal{Y}$. Although requiring significant 
calculation, \cref{alg:proxy} is only a pre-processing step and only needs to be applied once for one domain pair 
$\mathcal{X} \times \mathcal{Y}$. As described in the context of \cref{problem:targetmatrix}, in ID-based 
$\mathcal{H}^2$ construction, all the cluster pairs $\{X_i\times Y_i\}$ that are associated with sub-domains at the 
same level can fit to one domain pair $\mathcal{X} \times \mathcal{Y}$ by translations. Thus, 
for the compression of all these blocks $K(X_i, Y_i)$, the proxy point set $Y_p$ can be reused. 

%In addition, it is worth noting that the proxy-surface method described in \cref{sec:surface} is equivalent to 
%the proposed algorithm with $Y_p$ heuristically selected on the interior boundary $\Gamma$ of $\mathcal{Y}$. 
%We refer to this selection as {\it Surface Selection}. Just like \textit{Random Selection}, the size of $Y_p$ 
%needs to be manually decided in this selection method. 

%To summarize the algorithm, we first select a proxy point set $Y_p$ by any of the above methods. 
%With the target point set pair $X_0\times Y_0$, find an ID approximation of $K(X_0, Y_p)$
%with error threshold $\varepsilon$ and obtain the point subset $X_\text{rep} \subset X_0$ and matrix 
%$W_\text{rep}$. 
%An ID approximation of $K(X_0, Y_0)$ is then defined as in \cref{eqn:lowrankapprx} with 
%the approximation error for each row bounded as in \cref{ineqn:errorbound2} or \cref{ineqn:rowerrorbound}. 

\section{Comparison with existing methods}
In this section, we qualitatively compare our proposed algorithm with existing methods for 
the low-rank approximation of $K(X_0, Y_0)$ with $X_0 \times Y_0$ in some domain pair 
$\mathcal{X}\times \mathcal{Y}$.

In ID using sRRQR, each $K(x_i, y)$ is well approximated by $w_i^T K (X_\text{rep}, y)$ for $y\in Y_0$. 
The proposed algorithm, on the other hand, approximates each $K(x_i,y)$ by 
$w_i^T K (X_{\text{rep}}, y)$ for $y$ in the domain $\mathcal{Y}$ that contains $Y_0$.  
Thus, the proposed algorithm generally selects a larger row subset and the obtained 
$w_i$ does not necessarily minimize $\|K(x_i, Y_0) - w_i^T K (X_{\text{rep}}, Y_0)\|_2$.

%In addition, it is worth noting that the proxy-surface method described in \cref{sec:surface} is equivalent to 
%the proposed algorithm with $Y_p$ heuristically selected on the interior boundary $\Gamma$ of $\mathcal{Y}$. 
%We refer to this selection as {\it Surface Selection}. Just like \textit{Random Selection}, the size of $Y_p$ 
%needs to be manually decided in this selection method. 

It is also worth noting that the proxy-surface method described in \cref{sec:surface} is
equivalent to the proposed algorithm when $Y_p$ is chosen to discretize the interior boundary 
$\Gamma$ of $\mathcal{Y}$. We refer to this selection of $Y_p$ as {\it Surface Selection}.
Just like \textit{Random Selection}, the number of points in $Y_p$ needs to be manually decided 
in this selection scheme. For kernels from potential theory, it has been shown analytically in
\cref{sec:surface} that {\it Surface Selection} of $Y_p$ is sufficient for the proposed 
algorithm to be effective. 
However, for general kernels, this selection scheme usually leads to much larger error in comparison
with \textit{Random} and \textit{ID Selection} schemes in our numerical experiments.

\subsection{Comparison with algebraic and analytic methods} 
%For the following discussion, we will discuss the difference between the proposed
%algorithm and general algebraic and analytic methods in a united abstract framework.
%
In general, an analytic method seeks a degenerate approximation of $K (x, y)$ 
over $\mathcal{X} \times \mathcal{Y}$ as
\begin{equation}\label{eqn:expansion}
  K (x, y) \approx \sum_{i = 1}^r f_i (x) h_i (y), \quad x \in
  \mathcal{X}, y \in \mathcal{Y},
\end{equation}
with some analytically calculated functions $\{f_i(x)\}$ and $\{h_i(y)\}$. 
%The low-rank approximation is then obtained as $K (X_0, Y_0) \approx F (X_0)^T H (Y_0)$ 
%where $F (x) = (f_i (x))_{i = 1}^r$ and $H (y) = (h_i(y))_{i = 1}^y$ are both column vector
%functions. 
%For any $X_0 \times Y_0 \subset \mathcal{X} \times \mathcal{Y}$, the rank $r$ is fixed and
%the above approximation has bounded entry-wise error. 

An algebraic method, on the other hand, seeks a degenerate approximation of $K (x, y)$ over 
$X_0 \times Y_0$ instead. Denoting an obtained low-rank approximation as $K (X_0, Y_0) \approx F^T H$
with $F \in \mathbb{R}^{r \times | X_0 |}$ and $H \in \mathbb{R}^{r \times | Y_0 |}$, 
the degenerate approximation can be defined as \cref{eqn:expansion} with $\mathcal{X}\times\mathcal{Y}$ 
replaced by $X_0\times Y_0$ and 
\begin{equation}
  f_i (x) = \sum_{x_j \in X_0} F_{i j} \delta_{x_j} (x), \quad h_i (y) =
  \sum_{y_j \in Y_0} H_{i j} \delta_{y_j} (y),
\end{equation}
where $\delta_v (z) = 
\left\{ \begin{array}{ll}
  1, & \ z = v\\
  0, & \ z \neq v
\end{array} \right.$ is the delta function. 
%The optimal rank $r$, or equivalently the
%smallest number of expansion terms, is always smaller than that by analytical methods
%since the domain $X_0 \times Y_0$ is only a subset of $\mathcal{X} \times \mathcal{Y}$. 

Similarly, the proposed algorithm can be regarded as seeking a degenerate approximation of
$K (x, y)$ over $X_0 \times \mathcal{Y}$. The degenerate approximation 
can still be defined as \cref{eqn:expansion} with $\mathcal{X}\times\mathcal{Y}$ 
replaced by $X_0\times \mathcal{Y}$, the summation over $x_i \in X_\text{rep}$ and 
\begin{equation}
  f_i (x) = \sum_{x_j \in X_0} (W_\text{rep})_{ji} \delta_{x_j} (x), \quad h_i (y) = K
  (x_i, y).
\end{equation}
%where the summation is over $x_i \in X_\text{rep}$.
%The rank $r$ only depends on $X_0$ and can be adaptively determined by the prescribed error 
%threshold. 

Theoretically, the optimal degenerate approximation (i.e., the expansion with least 
terms to meet the same accuracy criteria) of $K(x,y)$ in a smaller domain pair should have 
fewer expansion terms. 
As $X_0\times Y_0 \subset X_0 \times \mathcal{Y} \subset \mathcal{X}\times \mathcal{Y}$, 
algebraic methods can obtain the smallest optimal rank while analytic methods deliver the 
largest optimal rank among the three classes of methods. 
Our proposed algorithm lies between analytic and algebraic methods and can be regarded as balancing
between computational cost, where analytic methods are better, and optimal rank estimation, 
where algebraic methods are better.

\subsection{Comparison with KIFMM}
Here, we focus on the source to multipole (S2M) translation phase in KIFMM. Similar comparisons with 
the other phases can be easily established.

\begin{figure}[ht]
	\centering
	\label{fig:illustration_kifmm}
	\includegraphics[width=0.8\textwidth]{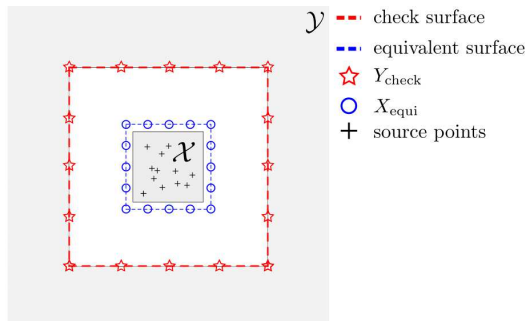}
	\caption{2D illustration of the S2M phase in KIFMM. The strong admissibility condition is applied over
	$\mathcal{X} \times \mathcal{Y}$. The equivalent surface lies between $\mathcal{X}$ and $\mathcal{Y}$
	and encloses $\mathcal{X}$, and the check surface lies between $\mathcal{Y}$ and the equivalent surface. 
	The potential at $y\in\mathcal{Y}$ generated by source charges in $\mathcal{X}$ is expected
	to be also generated by equivalent charges distributed on the equivalent surface. 
	The equivalent charges can be determined by matching their induced potential on the check surface with that 
	induced by the source charges. $X_\text{equi}$ and $Y_\text{check}$ are the discretization points
	of the associated surfaces.
	}
\end{figure}

An illustration of the S2M phase in KIFMM is shown in \cref{fig:illustration_kifmm}. 
Taking the Laplace kernel $K(x,y)$ in 2D as an example, the potential $q (y)$ at 
$y \in \mathcal{Y}$ generated by a source point set $X_0$ with charges $\{f_i\}$ is 
calculated as
\begin{equation}
  q (y) = \sum_{x_i \in X_0} K (x_i, y) f_i = K (X_0, y)^T F,
\end{equation}
where $K (X_0, y) = (K (x_i, y))_{x_i\in X_0}$ and $F = (f_i)_{x_i\in X_0}$ are column vectors.

KIFMM calculates the equivalent charges $\tilde{F}$ at grid points $X_{\text{equi}}$
on the equivalent surface by matching the potentials at grid points $Y_{\text{check}}$ 
on a check surface generated by both $F$ and $\tilde{F}$ as
\begin{equation}
  K (X_{\text{equi}}, Y_{\text{check}})^T \tilde{F} = K (X_0,
  Y_{\text{check}})^T F.
\end{equation}
%A regularized system is usually solved in KIFMM but the result is similar.
The potential at $y$ is then approximated as
\begin{align*}
  q (y) 
  &= K (X_0, y)^T F \\
  &\approx K (X_{\text{equi}}, y)^T \tilde{F} = K
  (X_{\text{equi}}, y)^T (K (X_{\text{equi}}, Y_{\text{check}})^T)^{\dag} K
  (X_0, Y_{\text{check}})^T F.
\end{align*}
This approximation applies to any charge vector $F$ and source points $X_0\subset\mathcal{X}$. 
Thus, it is equivalent to approximating $K (x, y)$ over $\mathcal{X}\times \mathcal{Y}$ as
\begin{equation} \label{apprx:kifmm}
  K (x, y) \approx K (x, Y_{\text{check}}) K (X_{\text{equi}},
  Y_{\text{check}})^{\dag} K (X_{\text{equi}}, y).
\end{equation}

For the proposed algorithm, combining the equations \cref{eqn:representation} and \cref{eqn:SYp},
it holds that
\begin{equation}\label{apprx:proxy}
  K(x, y) = K (x, Y_p) K (X_p, Y_p)^{\dag} K (X_p, y). 
\end{equation}
For non-degenerate kernels, the above equation becomes an approximation that shares exactly the 
same form as \cref{apprx:kifmm}. 

From the above discussion, the S2M phase in KIFMM and the proposed algorithm are 
based on a similar degenerate approximation of $K (x, y)$ over $\mathcal{X} \times \mathcal{Y}$.
However, in the proposed algorithm, $X_p$ and $Y_p$ are free to be selected in the whole domain 
pair $\mathcal{X} \times \mathcal{Y}$ while $X_{\text{equi}}$ and $Y_{\text{check}}$ are restricted
to be on the pre-defined equivalent surface and check surface respectively. 
In addition, the proposed algorithm only implicitly depends on \cref{apprx:proxy} and also takes into account
$X_0$ to automatically determine the rank $r$ for a given error threshold.
For the S2M phase in KIFMM, the rank of the underlying approximation  \cref{apprx:kifmm} is fixed to be
 $| X_{\text{equi}} |$ and needs to be manually decided. 
It should be expected that for general kernels where no Green's Theorem exists, \cref{apprx:kifmm} might 
not be accurate due to the restriction of the locations of $X_\text{equi}$ and $Y_\text{check}$, just like 
the proxy-surface method.

\section{Numerical experiments}
The basic settings for all the numerical experiments are as follows.

\begin{itemize}
\item Two kernels are tested: $K_1(x,y) = 1/|x-y|$ and $K_2(x,y) = \sqrt{1+|x-y|^2}$. 
$K_1(x,y)$ is the Laplace kernel in 3D where the proxy-surface method works. 

\item The tested domains are selected as follows with dimension $d = 2,3$.  
	\begin{itemize}
	\item Far-apart domain pair: $\mathcal{X} = [- 1, 1]^d, \mathcal{Y} = [-9, 9]^d \backslash [- 3, 3]^d$,
	\item Nearby domain pair: $\mathcal{X} = [- 1, 1]^d, \mathcal{Y} = [-9, 9]^d \backslash [-1.1, 1.1]^d$.
	\end{itemize} 
\item Point sets $X_0$ and $Y_0$ are all uniformly and randomly distributed 
within $\mathcal{X}$ and $\mathcal{Y}$, respectively, unless otherwise specified.

\item The error threshold for the ID approximation of $K(X_0, Y_p)$ is set as $10^{-6}\sqrt{|Y_p|}$ so
that the average entry-wise approximation error of each row satisfies
\[ 
  \| e_i (Y_p) \|_2 / \sqrt{|Y_p|} \leqslant 10^{- 6}, \quad \forall x_i \in X_0.
\] 

\item The entry-bound parameter $C$ for sRRQR in both the ID approximation of $K(X_0, Y_p)$ and 
\cref{alg:proxy} is set to 2. 

\item Denote the proxy point sets obtained by {\it Random, ID} and {\it Surface Selection} 
schemes as $Y_p^\text{rand}$, $Y_p^\text{id}$ and $Y_p^\text{surf}$ respectively. 
\cref{alg:proxy} for the selection of $Y_p^\text{id}$ has an initial point set pair $X_d\times Y_d$ of 
size approximately $1500\times 15000$. 
\cref{tab:proxy_point} lists the sizes of the selected $Y_p^\text{id}$ and the Matlab 
runtime of \cref{alg:proxy} with different settings.
\cref{fig:proxy} shows the selected $Y_p^\text{id}$ for the two kernels for the far-apart 
domain pair in 2D. 
The number of points in $Y_p^\text{rand}$ is set to 2000 based on the results of $|Y_p^\text{id}|$ in 
\cref{tab:proxy_point}.

\begin{table}[ht]
\centering
\caption{Proxy point set size $|Y_p^\text{id}|$ and runtime of \cref{alg:proxy} with different 
settings (\#\slash sec.). \label{tab:proxy_point} }
\begin{tabular}{c|cc|cc}
\toprule
		& 	\multicolumn{2}{c|}{$K_1(x,y)$} & \multicolumn{2}{c}{$K_2(x,y)$} \\
		& 	2D & 3D & 2D & 3D \\
\midrule
Far-apart domain & 174/5.35 & 636/9.05 & 126/7.83 & 821/13.32\\
Nearby domain    & 500/4.98 & 797/8.93 & 316/10.90 & 991/14.16\\
\bottomrule
\end{tabular}
\end{table}

\begin{figure}[ht]
\centering
\subfloat[$K_1(x,y)\!=\!1/|x\!-\!y|$ with $|Y_p^\text{id}| = 174$]{
\includegraphics[width=0.43\textwidth]{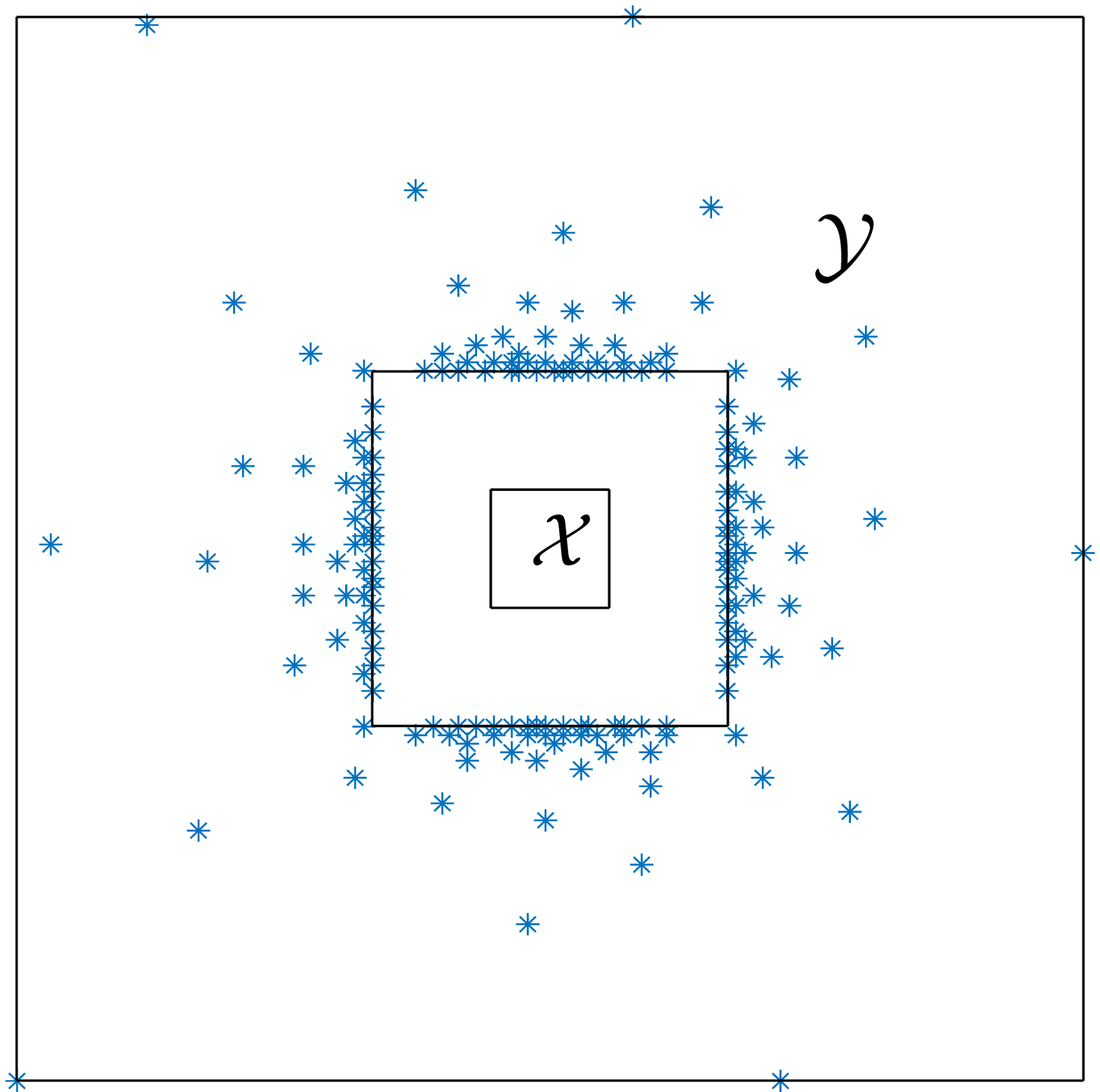}
}
\subfloat[$K_2(x,y)\!=\!\sqrt{1\!+\!|x\!-\!y|^2}$ with $|Y_p^\text{id}| = 126$]{
\includegraphics[width=0.43\textwidth]{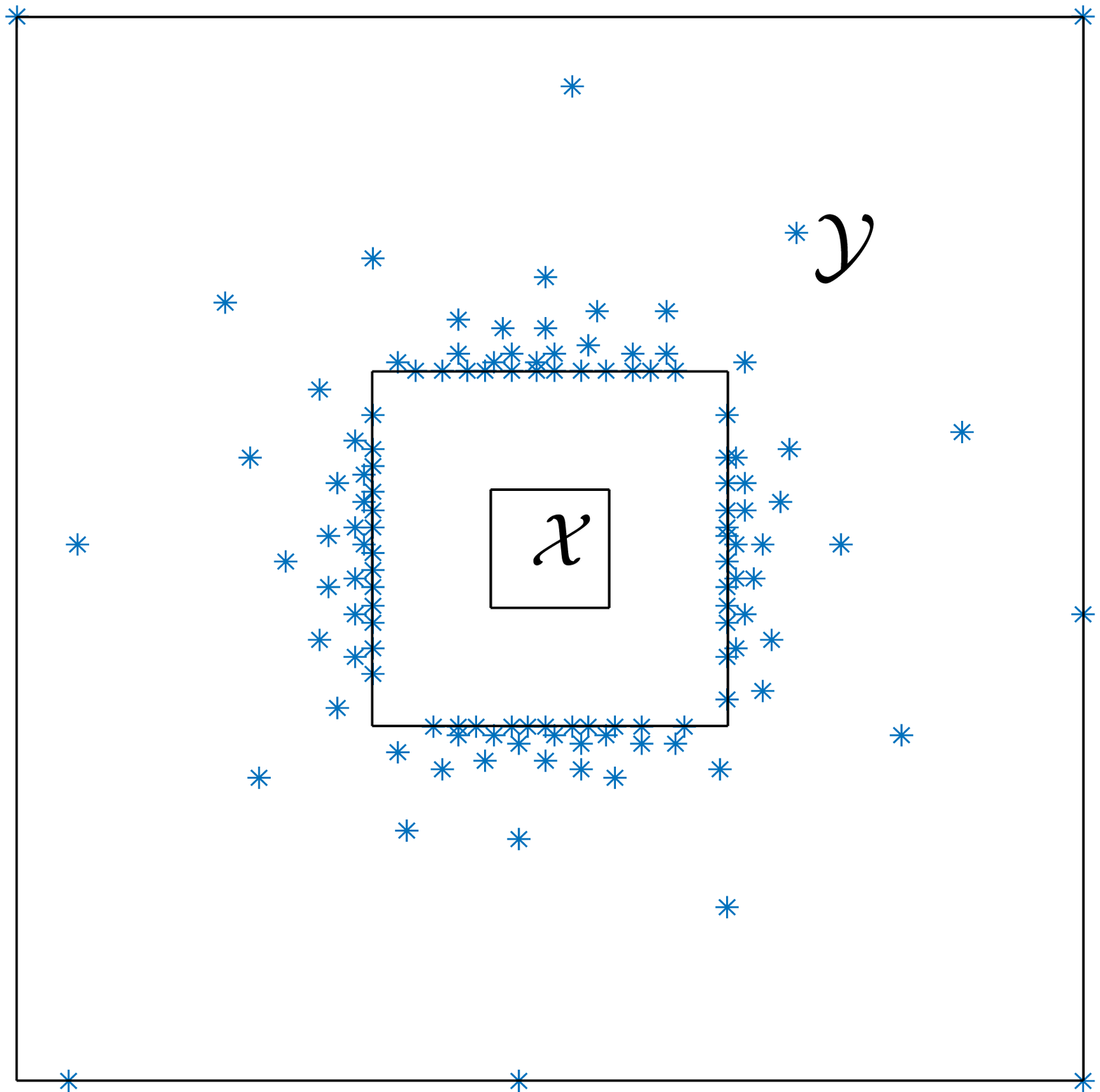}
}
\caption{Proxy point set $Y_p^\text{id}$ for two kernels with far-apart domain pair
$\mathcal{X}\times \mathcal{Y}$ in 2D. \label{fig:proxy}}
\end{figure}
\end{itemize}

%%%%%%%%%%%%%%%%%%%%%%%%%%%%%%%%%%%%%%%%%%%%%%%%%%%%%%%%%%%%%%%%%
%   Numerical Test 1
%%%%%%%%%%%%%%%%%%%%%%%%%%%%%%%%%%%%%%%%%%%%%%%%%%%%%%%%%%%%%%%%
\subsection{Function approximation error $|e_i(y)|$\label{sec:test1}}
The accuracy of the proposed approximation can be quantified by the function approximation error 
$e_i(y)$ in $\mathcal{Y}$ that is connected to $e_i(Y_p)$ as $e_i(y) = e_i(Y_p)S_{Y_p}(y)$ 
from \cref{eqn:error_representation}.
Viewing a general kernel as a degenerate kernel plus a small remainder,
the connection \cref{eqn:error_representation} only holds approximately.  
In this subsection, we check the obtained error $|e_i(y)|$ of applying the proposed algorithm
to the prescribed non-degenerate kernels. 

Consider the far-apart domain pair in 2D and $X_0\subset\mathcal{X}$ with 1000 points. 
With the prescribed error threshold and $Y_p^\text{id}$, the proposed algorithm obtains $X_\text{rep}$ 
with 35 points for $K_1(x,y)$ and 33 points for $K_2(x,y)$. Similarly, with $Y_p^\text{rand}$, 
the obtained $X_\text{rep}$ has 28 points for $K_1(x,y)$ and 29 points for $K_2(x,y)$. 
To check $e_i(y)$, a dense uniform grid in $\mathcal{Y}$ with 
approximately 40000 points is defined as $Y_0$.
 
By selecting a large point set $X_p$ in $\mathcal{X}$ with $|X_p| \gg |Y_p|$, $S_{Y_p} (Y_0)$ can be 
explicitly estimated as $K (X_p, Y_p)^{\dag} K (X_p, Y_0)$ by \cref{eqn:SYp} and the maximum errors of 
\cref{eqn:error_representation} for the two kernels are then found as
\[
\max_{x_i\in X_0, y\in Y_0} |e_i(y) - e_i(Y_p)S_{Y_p}(y)|
= \left\{
\begin{array}{ll}
4.20\times 10^{-11} & \text{ for } K_1(x,y) \text{ with } Y_p^\text{id}  \\
5.04\times 10^{-10} & \text{ for } K_1(x,y) \text{ with } Y_p^\text{rand} \\
2.39\times 10^{-10} & \text{ for } K_2(x,y) \text{ with } Y_p^\text{id} \\
7.30\times 10^{-9}  & \text{ for } K_2(x,y) \text{ with } Y_p^\text{rand} \\
\end{array}
\right.
.
\]
By these results, the connection \cref{eqn:error_representation} indeed holds approximately and 
thus the proposed algorithm should work with these non-degenerate kernels. 

\begin{figure}[ht]
	\centering
	\vspace{-0.4em}
	\includegraphics[width=\textwidth]{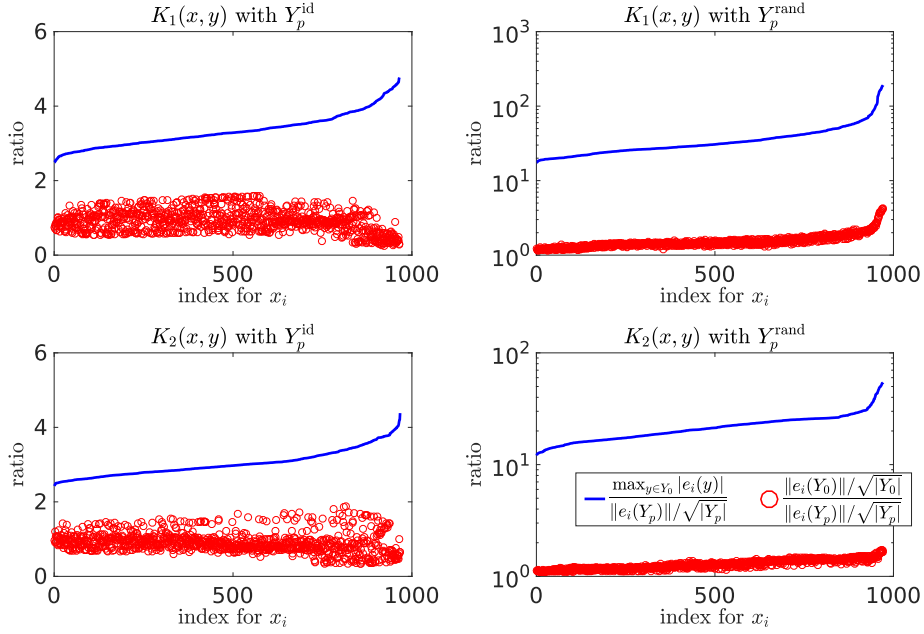}
	\vspace{-1.7em}
	\caption{Entry-wise ratios \cref{def:ratio} for two kernels with the far-apart domain  
	pair in 2D. Indices for $x_i \in X_0 \backslash X_\text{rep}$ are sorted such that
	the maximum ratios are in ascending order.\label{fig:test1_ratio}}
\end{figure}

To check the approximation $K(x_i, y) \approx w_i^T K(X_\text{rep}, y)$ for $x_i \in 
X_0 \backslash X_\text{rep}$\footnotemark, 
\footnotetext{The function approximation is exact for any $x_i \in X_\text{rep}$.}
\cref{fig:test1_ratio} shows the following entry-wise error ratios with both $Y_p^\text{id}$ and $Y_p^\text{rand}$,
\begin{equation}\label{def:ratio}
 	\dfrac{\max_{y\in Y_0}|e_i(y)|}{\|e_i(Y_p)\|_2/\sqrt{|Y_p|}} \quad \text{and}\quad 
 	\dfrac{\|e_i(Y_0)\|_2/\sqrt{|Y_0|}}{\|e_i(Y_p)\|_2/\sqrt{|Y_p|}}, \quad
 	\text{for any } x_i \in X_0 \backslash X_\text{rep}.
\end{equation}
%\begin{figure}[ht]
%	\centering
%	\subfloat[$K_1(x,y) = 1/ |x-y|$]
%	{\includegraphics[width=0.48\textwidth]{./figure/test1_ratio_strong_inverse}}
%	\subfloat[$K_2(x,y) = \sqrt{1+|x-y|^2}$]
%	{\includegraphics[width=0.48\textwidth]{./figure/test1_ratio_strong_mq}}
%	\caption{Entry-wise ratio as in \cref{def:ratio} for two kernels with the far-apart domain setting in 2D. 
%	Indices for $x_i \in X_0 \backslash X_\text{rep}$ are sorted to make 
%		the maximum ratios in ascending order.\label{fig:test1_ratio}}
%\end{figure}

By the previous analysis, the entry-wise error ratios in \cref{def:ratio} can be bounded as
\[
\dfrac{\|e_i(Y_0)\|_2/\sqrt{|Y_0|}}{\|e_i(Y_p)\|_2/\sqrt{|Y_p|}} \leqslant
\dfrac{\max_{y\in Y_0}|e_i(y)|}{\|e_i(Y_p)\|_2/\sqrt{|Y_p|}} \leqslant
\sqrt{|Y_p|}\max_{y\in Y_0}\|S_{Y_p}(y)\|_2 \leqslant
|Y_p|C
\]
where the last inequality only holds for $Y_p^\text{id}$. 
From the results in \cref{fig:test1_ratio}, the approximation $K(x_i,y)\approx w_i^TK(X_\text{rep},y)$ obtained 
by $Y_p^\text{id}$ for any $x_i \in X_0\backslash X_\text{rep}$ has its maximum error $e_i(y)$ in $\mathcal{Y}$ of 
the same scale as $e_i(Y_p)/\sqrt{|Y_p|}$ which is bounded by $10^{-6}$.
Thus, given any $Y_0\subset \mathcal{Y}$, the proposed approximation of $K(X_0, Y_0)$ should have entry-wise error
of the same scale as $10^{-6}$. 
Also, these results suggest that the above analytic upper bound may not be sharp for $Y_p^\text{id}$. 

From the results for $Y_p^\text{rand}$, $\|S_{Y_p^\text{rand}}(y)\|_2$ is large for some  
$y\in \mathcal{Y}$ that leads to very large $e_i(y)$. However, the average entry-wise error is still 
of scale $10^{-6}$ for each $e_i(y)$, which indicates that $\|S_{Y_p^\text{rand}}(y)\|_2$ may be small for $y$
in most of the domain $\mathcal{Y}$. Lastly, it is worth noting that \textit{ID Selection} obtains much better 
results with much fewer proxy points compared to \textit{Random Selection}.

%%%%%%%%%%%%%%%%%%%%%%%%%%%%%%%%%%%%%%%%%%%%%%%%%%%%%%%%%%%%%%%%%%%%
%   Test
%%%%%%%%%%%%%%%%%%%%%%%%%%%%%%%%%%%%%%%%%%%%%%%%%%%%%%%%%%%%%%%%%%%%
\subsection{Comparison with algebraic methods}
With a fixed cluster set $X_0$ and the prescribed error threshold, assume that the proposed algorithm
with $Y_p^\text{id}$ gives a rank-$|X_\text{rep}|$ approximation 
$K (X_0, Y_0) \approx W_\text{rep} K (X_\text{rep}, Y_0)$ for any $Y_0 \subset \mathcal{Y}$. 
We compare this approximation with those of the following methods:
\begin{itemize}
\item {\it The proposed algorithm with $Y_p^\textnormal{rand}$ and fixed rank $|X_\text{rep}|$}
\item {\it ID with row subset $X_\text{rep}$},
\item {\it ID using sRRQR with fixed rank $|X_\text{rep}|$},
\item {\it SVD wtih fixed rank $|X_\text{rep}|$},
\item {\it ACA with fixed rank $|X_\text{rep}|$}.
\end{itemize}
The proposed algorithm with a fixed rank $r$ means to find a rank-$r$ ID approximation of $K(X_0, Y_p)$ 
in the first step of the algorithm.  
ID with row subset $X_\text{rep}$ simply replaces $W_\text{rep}$ by $K(X_0, Y_0)K(X_\text{rep}, Y_0)^\dag$. 
ACA is implemented with partial pivoting as described in \cite{bebendorf_adaptive_2003}. 
%The proxy-surface method can be regarded as the proposed algorithm with the proxy point set 
%$Y_p$ selected over the interior surface $\Gamma$ of $\mathcal{Y}$. 

\begin{figure}[h]
    \centering
    \vspace{-1em}
    \subfloat[$K_1(x,y) = 1/|x-y|$]
    {\includegraphics[width=0.45\textwidth]{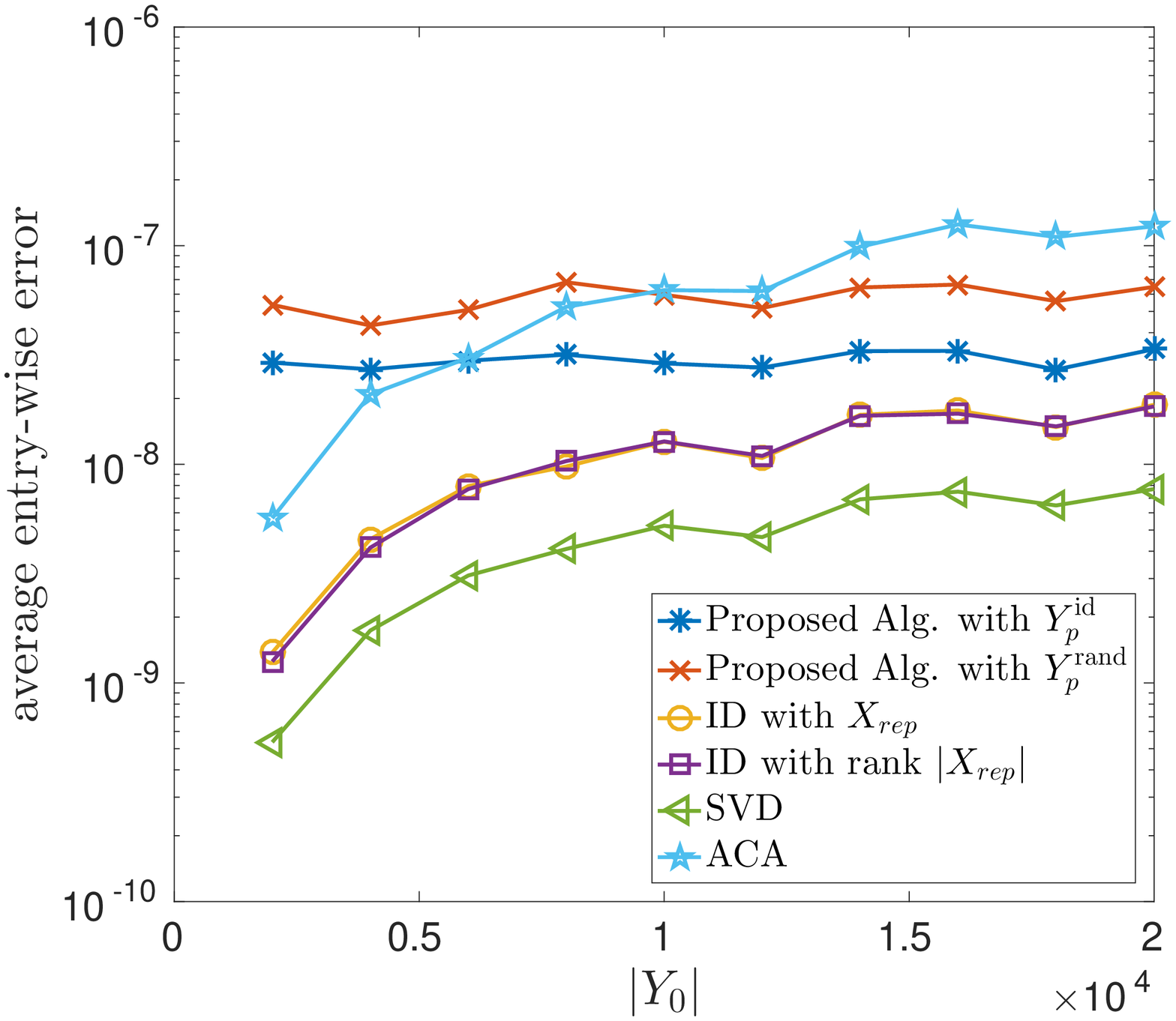}}
    \subfloat[$K_2(x,y) = \sqrt{1+|x-y|^2}$]
    {\includegraphics[width=0.45\textwidth]{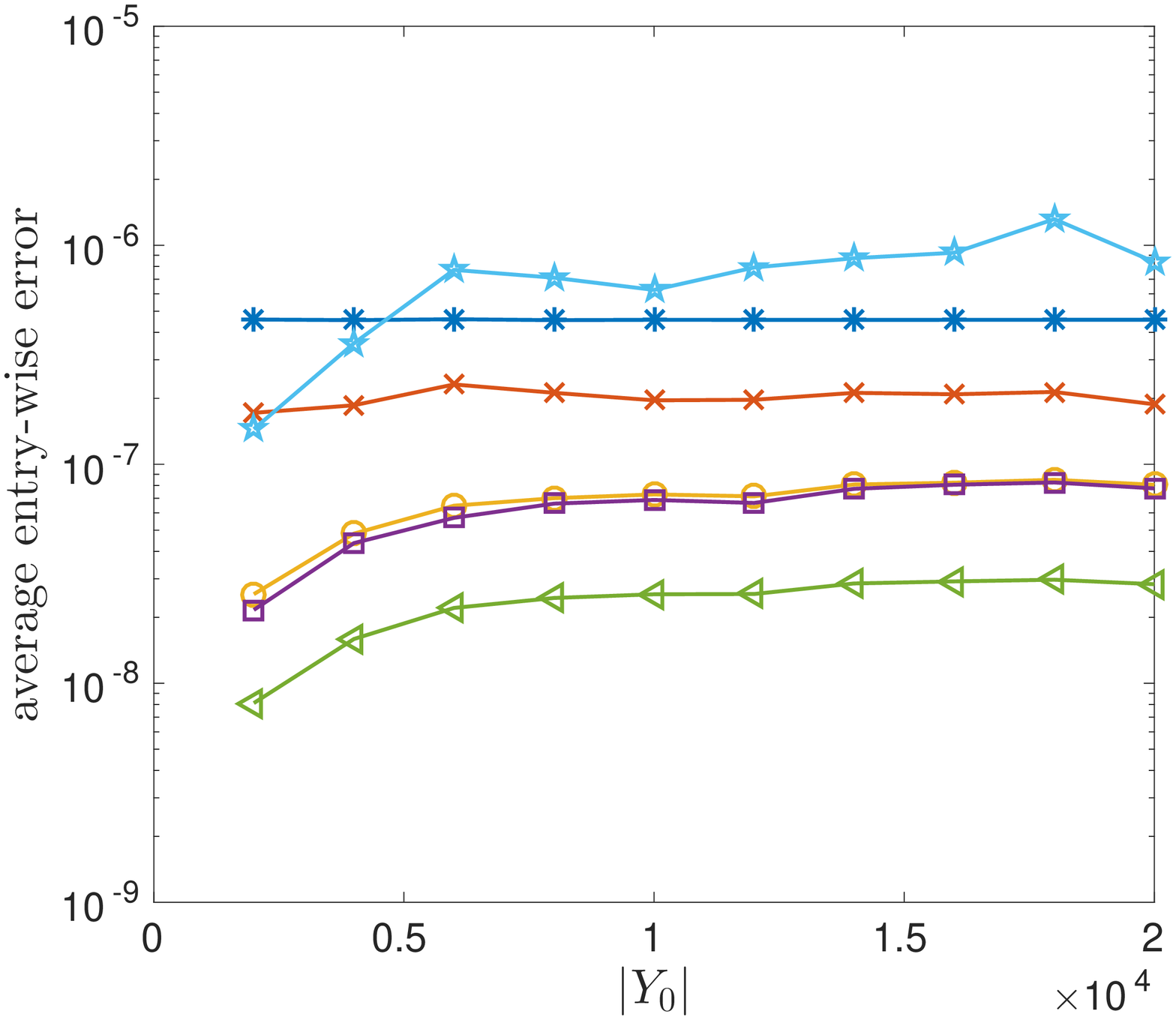}}
    \caption{Average entry-wise error of the obtained low-rank approximations.\label{fig:test2_error}}
\end{figure}

\begin{figure}[h]
    \centering
    \vspace{-1em}
    \subfloat[$K_1(x,y) = 1/|x-y|$]
    {\includegraphics[width=0.45\textwidth]{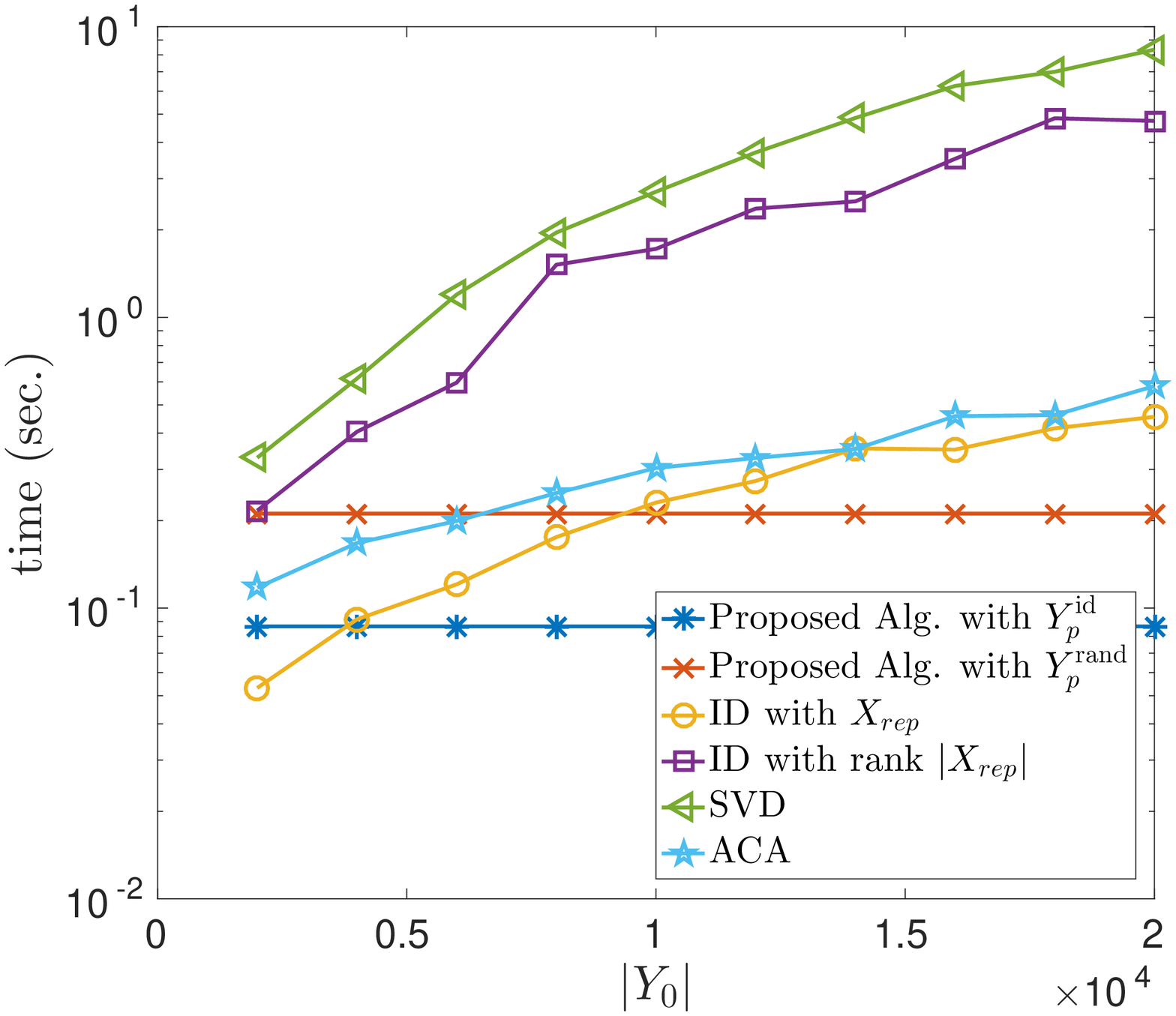}}
    \subfloat[$K_2(x,y) = \sqrt{1+|x-y|^2}$]
    {\includegraphics[width=0.45\textwidth]{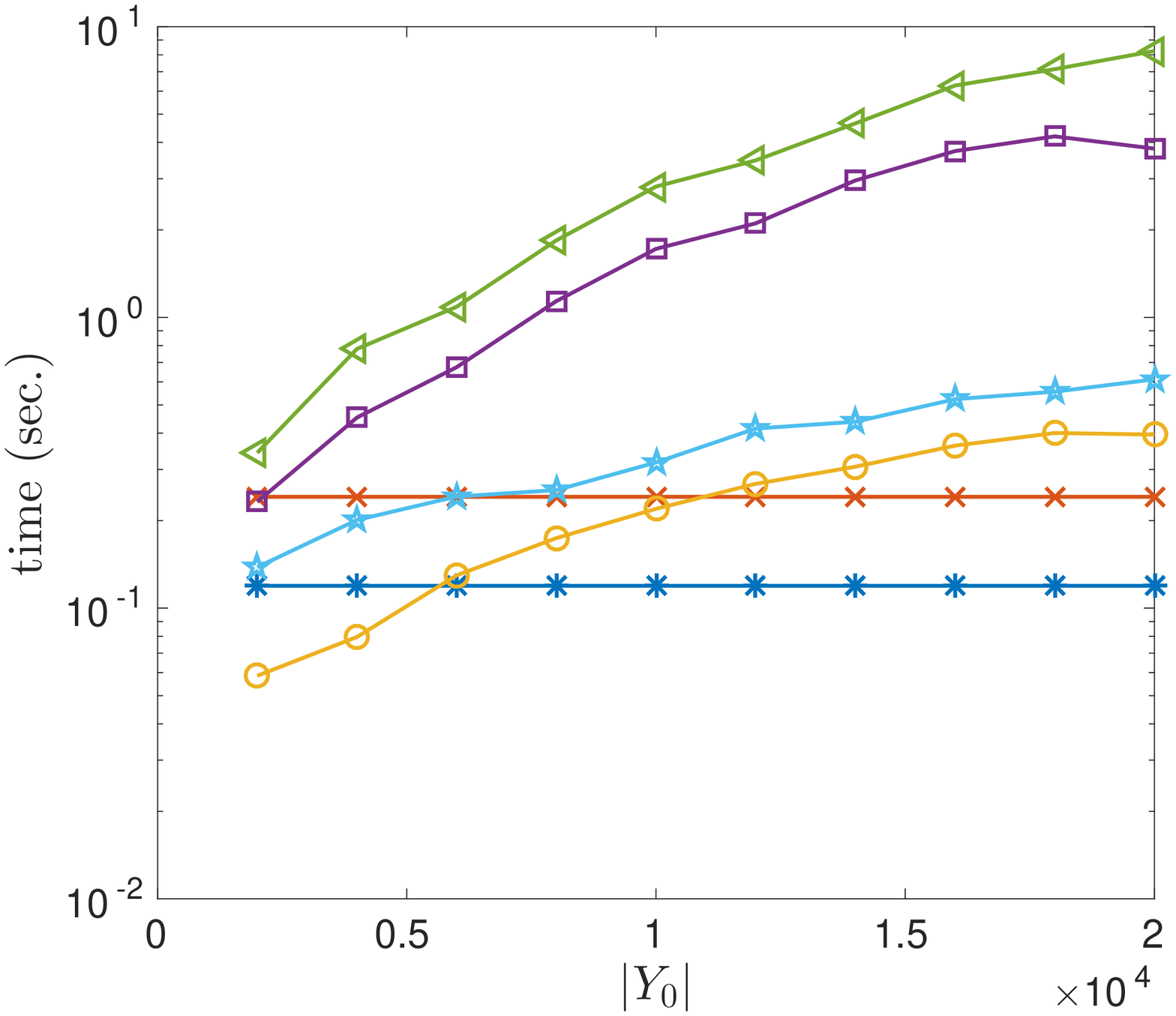}}
    \caption{Runtime of different low-rank approximation methods. The selection of $Y_p^\text{id}$ 
    by \cref{alg:proxy} takes $9.05$ sec.\ for $K_1(x,y)$ and $13.32$ sec.\ for $K_2(x,y)$.
    \label{fig:test2_time}}
\end{figure}

Consider the far-apart domain pair in 3D and $X_0 \subset \mathcal{X}$ with 1000 points. 
The obtained $X_\text{rep}$ has 119 points for $K_1(x,y)$ and 131 points for $K_2(x,y)$. 
Selecting $Y_0$ in $\mathcal{Y}$ with different number of points, 
\cref{fig:test2_error} shows the average entry-wise error of the low-rank approximations and
\cref{fig:test2_time} shows the runtime of our Matlab implementation. 


For the two kernels, the average entry-wise errors of the proposed algorithm with $Y_p^\text{id}$ and 
$Y_p^\text{rand}$ both remain roughly constant for different $|Y_0|$ and are close to those of the ID 
approximation using sRRQR.
The runtime of the proposed algorithm is independent of $|Y_0|$ which becomes advantageous 
over purely algebraic methods when $|Y_0|$ is large.
 
%It is interesting to note that the proxy-surface method also has a roughly constant-scale 
%entry-wise error which should be expected since it is our proposed algorithm with a specific 
%choice of $Y_p$. For $K_1(x,y) = 1/|x-y|$ from potential theory, 
%proxy-surface method has slightly smaller approximation error than the proposed algorithm, 
%indicating that the selection of $Y_p$ over the surface $\Gamma$ of $\mathcal{Y}$ is better 
%than that of \cref{alg:proxy} in this case. However, the proxy-surface method needs to 
%heuristically decide the size of $|Y_p|$ and insufficient points have been noted to 
%cause larger approximation error in our experiments. In addition, there is also an obvious
%accuracy deterioration of the proxy-surface method when applied over $K_2(x,y)$ that is not from
%potential theory. 

%%%%%%%%%%%%%%%%%%%%%%%%%%%%%%%%%%%%%%%%%%%%%%%%%%%%%%%%%%%%%%%%%%%%
%   Test
%%%%%%%%%%%%%%%%%%%%%%%%%%%%%%%%%%%%%%%%%%%%%%%%%%%%%%%%%%%%%%%%%%%%
%
\subsection{Comparison with different selections of $Y_p$\label{sec:compare_selection}}
The number of points in both $Y_p^\text{rand}$ and $Y_p^\text{surf}$ need to be manually specified 
while \cref{alg:proxy} determines the number of points in $Y_p^\text{id}$ automatically for different
kernels and domain pairs. 
Intuitively, $|Y_p^\text{id}|$, as an estimate of $r_{KL}$, should be a good reference value for 
$|Y_p^\text{rand}|$ and $|Y_p^\text{surf}|$. 
Assuming that the proposed algorithm with $Y_p^\text{id}$ and the prescribed error threshold gives 
$K(X_0, Y_0) \approx W_\text{rep} K(X_\text{rep}, Y_0)$, 
we compare the rank-$|X_\text{rep}|$ ID approximations obtained by the proposed algorithm with  
following proxy point sets, 
\begin{itemize}
\item $Y_p^\text{rand}$ with $\frac{1}{2}|Y_p^\text{id}|$, $|Y_p^\text{id}|$ and $2|Y_p^\text{id}|$ points.
\item $Y_p^\text{surf}$ with $\frac{1}{2}|Y_p^\text{id}|$, $|Y_p^\text{id}|$ and $2|Y_p^\text{id}|$ points.
\end{itemize}

\begin{figure}[h]
    \centering
    \vspace{-1em}
    \subfloat[$K_1(x,y) = 1/|x-y|$]
    {\includegraphics[width=0.45\textwidth]{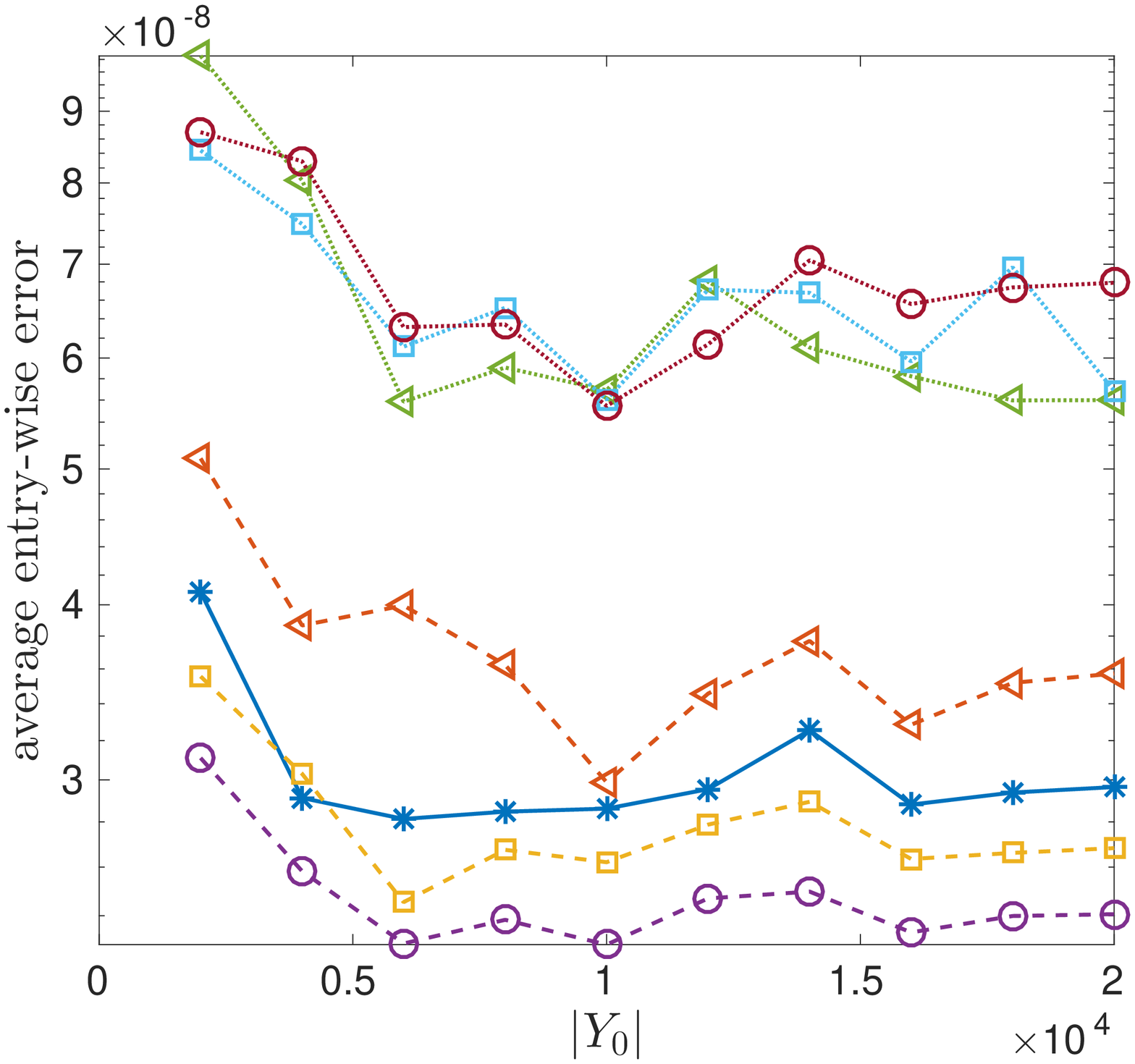}}
    \subfloat[$K_2(x,y) = \sqrt{1+|x-y|^2}$\label{fig:proxy_3d_strong_b}]
    {\includegraphics[width=0.45\textwidth]{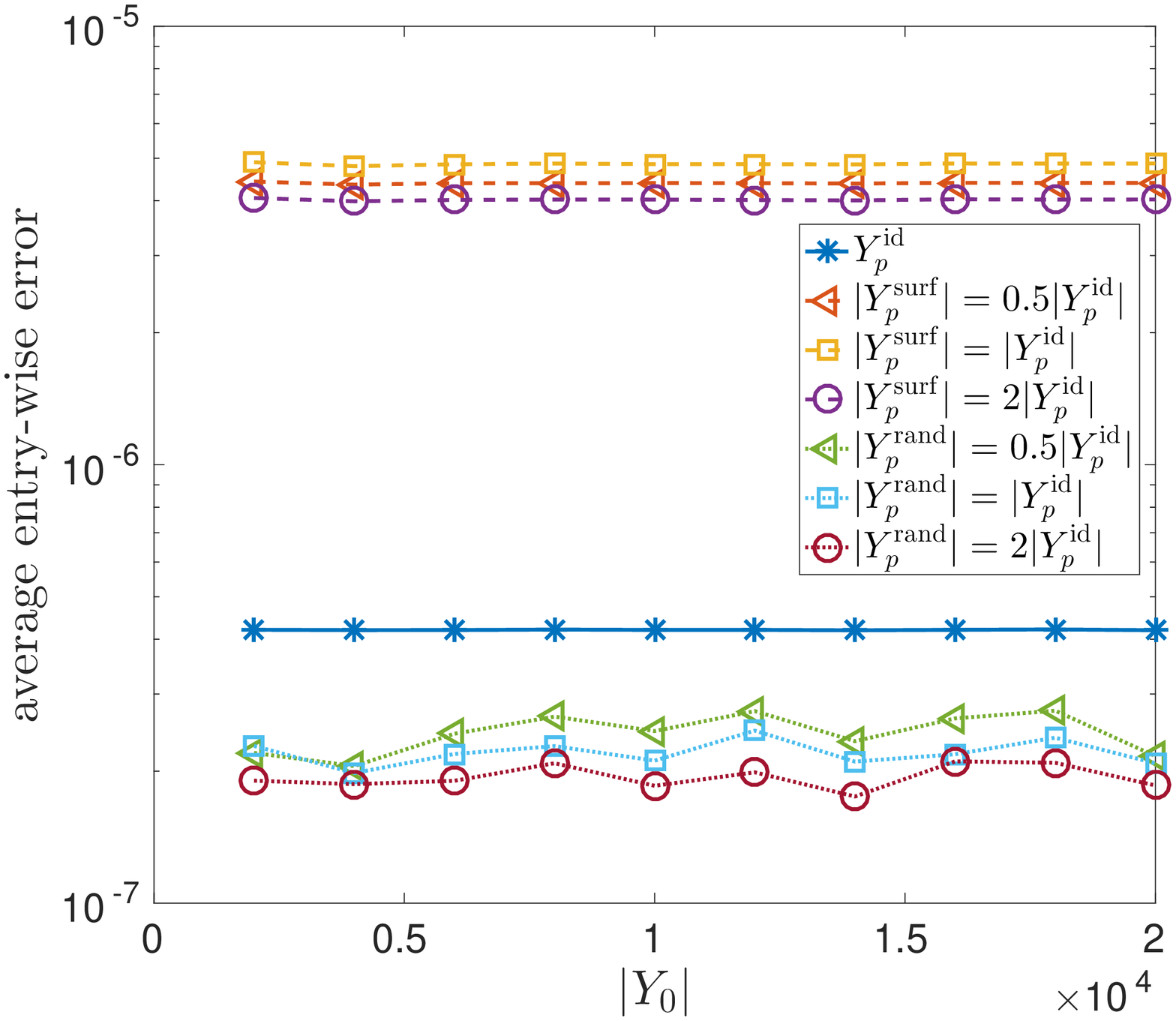}}
    \caption{Average entry-wise approximation error for two kernels with the far-apart domain 
    pair in 3D. For $K_1(x,y)$, $|Y_p^\text{id}| = 636$ and $|X_\text{rep}| = 135$. For $K_2(x,y)$, 
            $|Y_p^\text{id}|=821$ and $|X_\text{rep}|= 142$.\label{fig:proxy_3d_strong}}
\end{figure}

We continue to consider the far-apart domain pair in 3D and $X_0\in\mathcal{X}$ with 1000 points. 
\cref{fig:proxy_3d_strong} shows the average entry-wise approximation error with different numbers of points in $Y_0$. 
The proxy-surface method (i.e., $Y_p^\text{surf}$) gives the best approximations for the 3D Laplace kernel $K_1(x,y)$
while its accuracy degrades drastically for the general kernel $K_2(x,y)$. 
Thus, the analytic method here leads to a better selection of the proxy points but the method is only limited to 
specific kernels and may be counter-productive otherwise.  

\textit{Random Selection} gives similar errors to that of \textit{ID Selection} for both kernels, which can 
also be observed in \cref{fig:test2_error} of the previous test. 
This observation suggests that in some cases, \textit{Random Selection} can be a better alternative
of $Y_p^\text{id}$ since it requires no significant pre-calculation and it can be adapted to different
domain pairs $\mathcal{X}\times \mathcal{Y}$ easily.
However, it remains to decide a proper number of the points in $Y_p^\text{rand}$ since an insufficient number
of points can lead to larger errors as can be somewhat suggested by \cref{fig:proxy_3d_strong_b} and 
an excessive number of points can lead to more computation. Furthermore, we remind readers of the results in 
\cref{sec:test1} that at some $y \in \mathcal{Y}$, the entry-wise error $e_i(y)$ can be 10 or more
times larger than the expected error threshold $10^{-6}$. Thus, \textit{Random Selection} can have much 
worse performance than \textit{ID Selection} for $Y_0$ with specific point distributions.

\begin{figure}[h]
    \centering
    \vspace{-1em}
    \subfloat[$K_1(x,y) = 1/|x-y|$\label{fig:proxy_2d_weak_a}]
    {\includegraphics[width=0.45\textwidth]{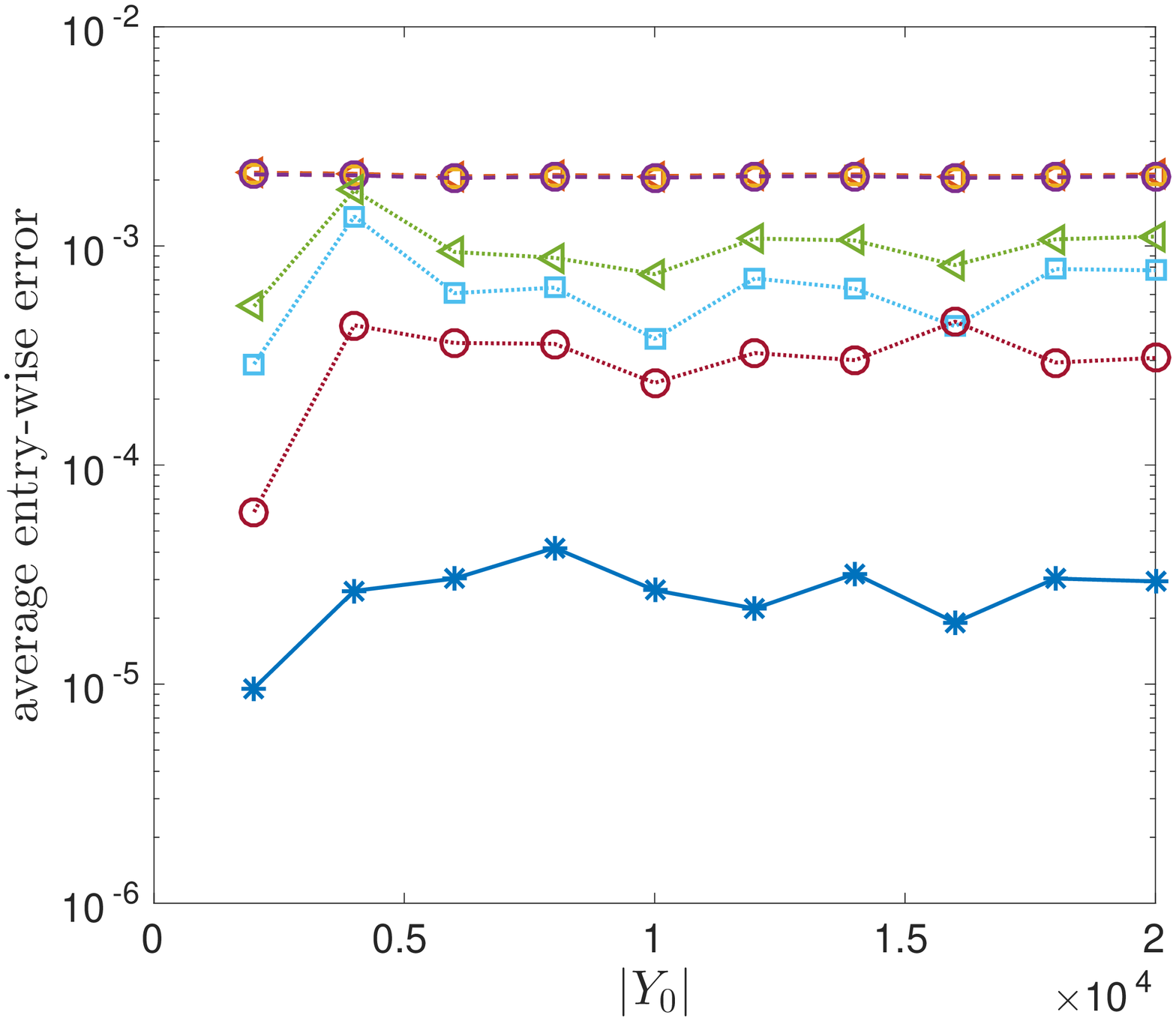}}
    \subfloat[$K_2(x,y) = \sqrt{1+|x-y|^2}$]
    {\includegraphics[width=0.45\textwidth]{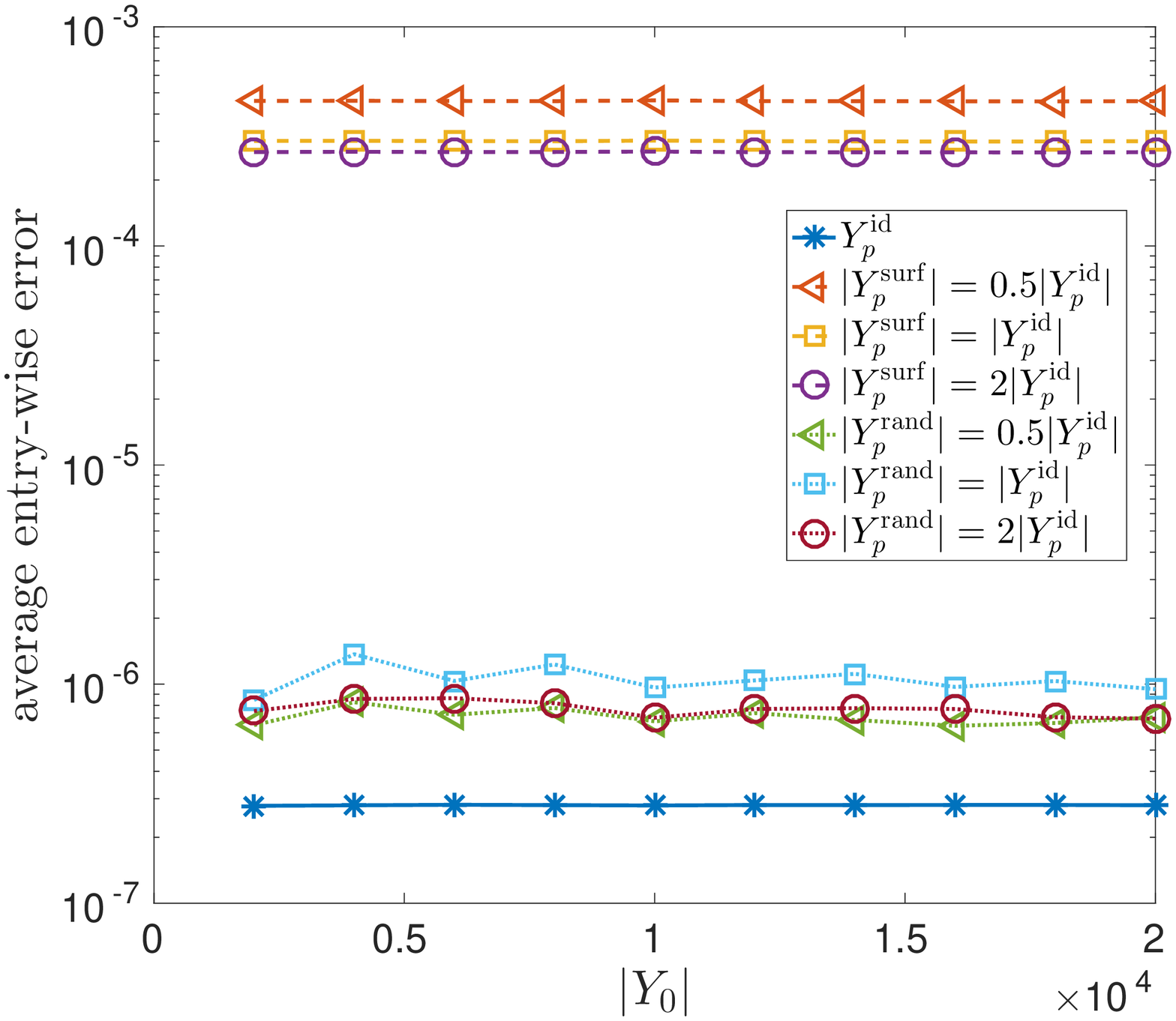}}
    \caption{Average entry-wise approximation error for two kernels with the nearby domain pair 
        in 2D. For $K_1(x,y)$, $|Y_p^\text{id}| = 500$ and $|X_\text{rep}| = 197$. For $K_2(x,y)$, 
        $|Y_p^\text{id}|=316$ and $|X_\text{rep}|= 72$. \label{fig:proxy_2d_weak}}
\end{figure}

More distinguishable differences between results from \textit{Random, ID} and \textit{Surface Selection} schemes
can be found for the two kernels with the nearby domain pair in 2D as shown in \cref{fig:proxy_2d_weak}. 
In these results, it should be noted that for $K_1(x,y)$, all the $Y_p$ selections give much larger error than 
$10^{-6}$ while the rank-$|X_\text{rep}|$ truncated SVD for any of the tested $Y_0$ has average entry-wise error
at most $10^{-7}$. The main cause of this accuracy degradation is the singularity of $K_1(x,y)$ when $x$ and $y$
are close. The analysis and solution for this problem will be explained in the next subsection.

%%%%%%%%%%%%%%%%%%%%%%%%%%%%%%%%%%%%%%%%%%%%%%%%%%%%%%%%%%%%%%%%%
%   Numerical Test
%%%%%%%%%%%%%%%%%%%%%%%%%%%%%%%%%%%%%%%%%%%%%%%%%%%%%%%%%%%%%%%%%
\subsection{Improvement of the proxy point selection}\label{sec:improvement}
\begin{figure}[h]
\centering
\subfloat[Proxy point set $Y_p^\text{id}$ with 441 points in ${[-2,2]^2}$ and 59 points outside.]
{\includegraphics[height=0.23\textheight]{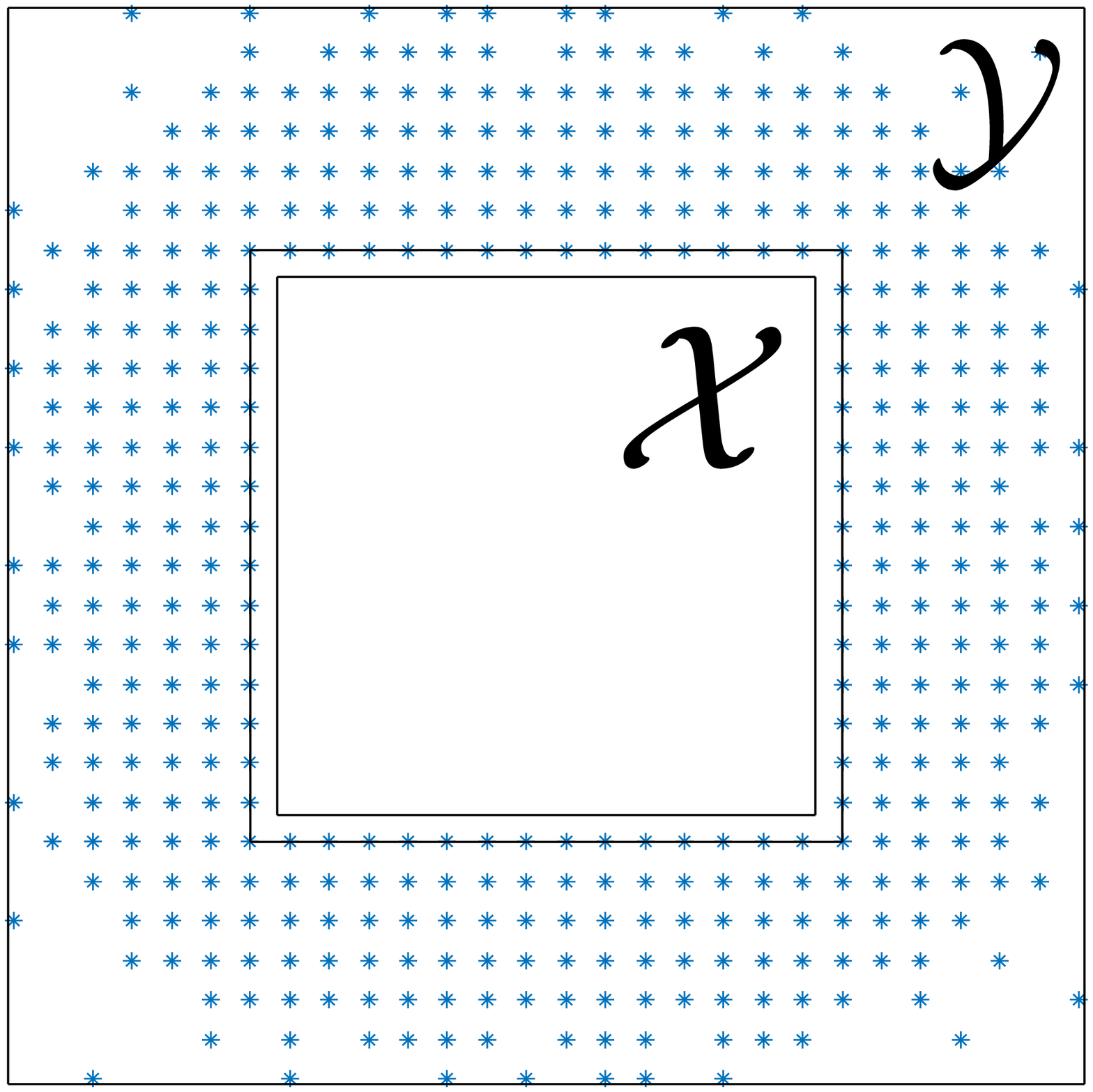}}
\hspace{1em}
\subfloat[$\frac{1}{|X_0|-|X_\text{rep}|}\sum_{x_i\in X_0\backslash X_\text{rep}}|e_i(y)|$ for ${y \in \mathcal{Y}}$ in $\log_{10}$ scale.]
{\includegraphics[height=0.23\textheight]{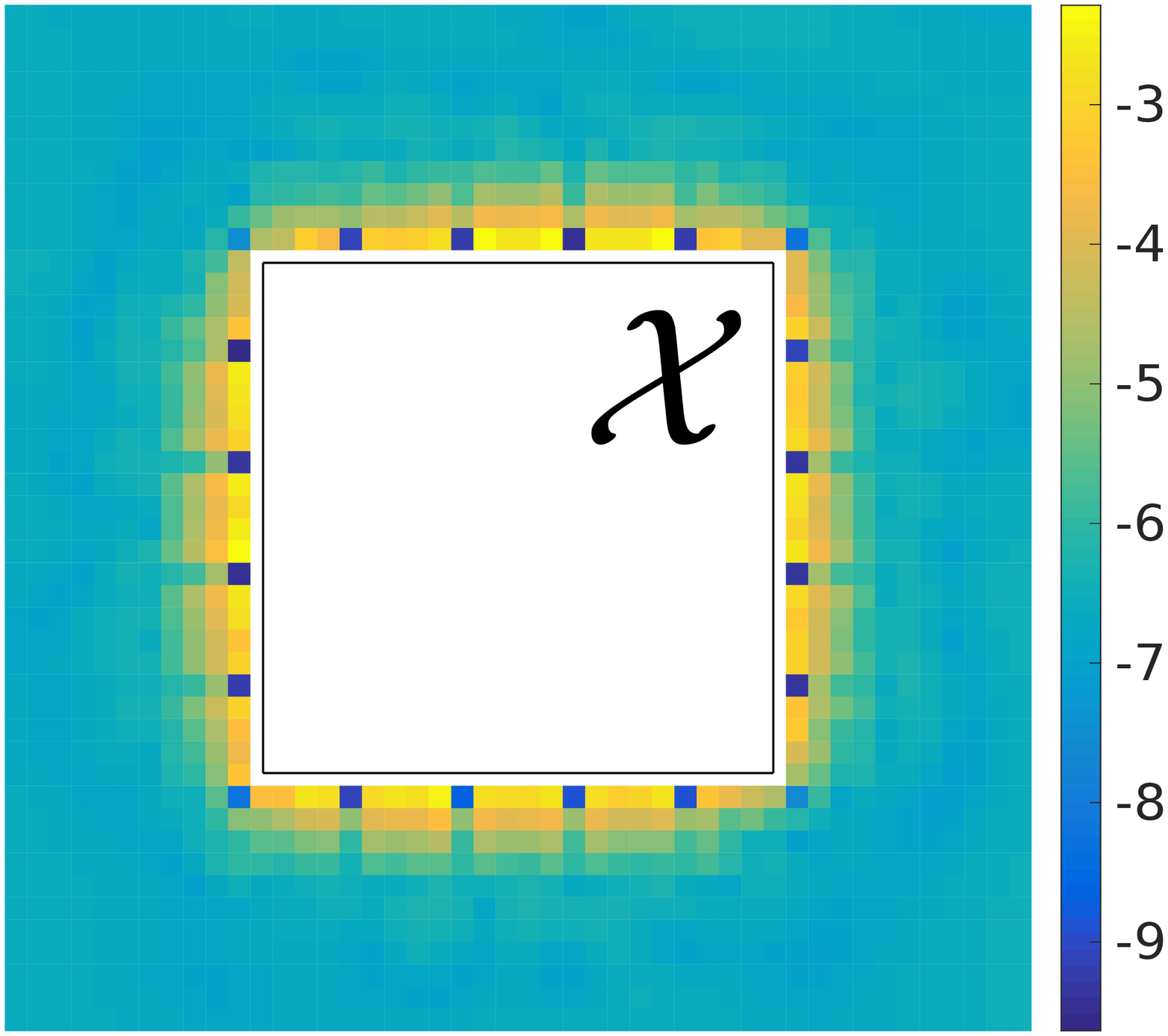}}
\caption{$Y_p^\text{id}$ obtained by \cref{alg:proxy} and the associated error distribution in $\mathcal{Y}$, 
	for $K_1(x,y) = \frac{1}{|x-y|}$ with the nearby domain pair in 2D.
The plotted area is $[-2,2]^2$ to highlight differences. The error function values outside of $[-2,2]^2$ are 
approximately $10^{-7}$. \label{fig:test4_weak}}
\end{figure}

Here, we introduce ideas to improve \textit{ID Selection}. The ideas can be easily adapted to 
improve \textit{Random Selection}. 
To better understand the large errors in \cref{fig:proxy_2d_weak_a}, the average error  
function $\frac{1}{|X_0|-|X_\text{rep}|}\sum_{x_i\in X_0\backslash X_\text{rep}}|e_i(y)|$ over $\mathcal{Y}$
and the selected proxy point set $Y_p^\text{id}$ are drawn in \cref{fig:test4_weak}. 

As can be observed, the largest errors are only located in the part of the domain near $\mathcal{X}$.
The most likely cause is that $\frac{1}{|x-y|}$ varies rapidly when $x$ and $y$ become close,
the candidate point set $Y_d$ in \cref{alg:proxy} may not be dense enough in the area near 
$\mathcal{X}$ to satisfy the prerequisite that $\text{col}(\Phi(Y_d)) = \mathbb{R}^{r_{KL}}$ in \cref{prop:Ysample}. 
A hint towards this is that most of the candidate points near $\mathcal{X}$ are selected for $Y_p^\text{id}$. 

\begin{figure}[h]
	\centering
	\subfloat[Proxy point set $Y_p^\text{id}$ with 879 points in ${[-2,2]^2}$ and 30 points outside.]
	{\includegraphics[height=0.23\textheight]{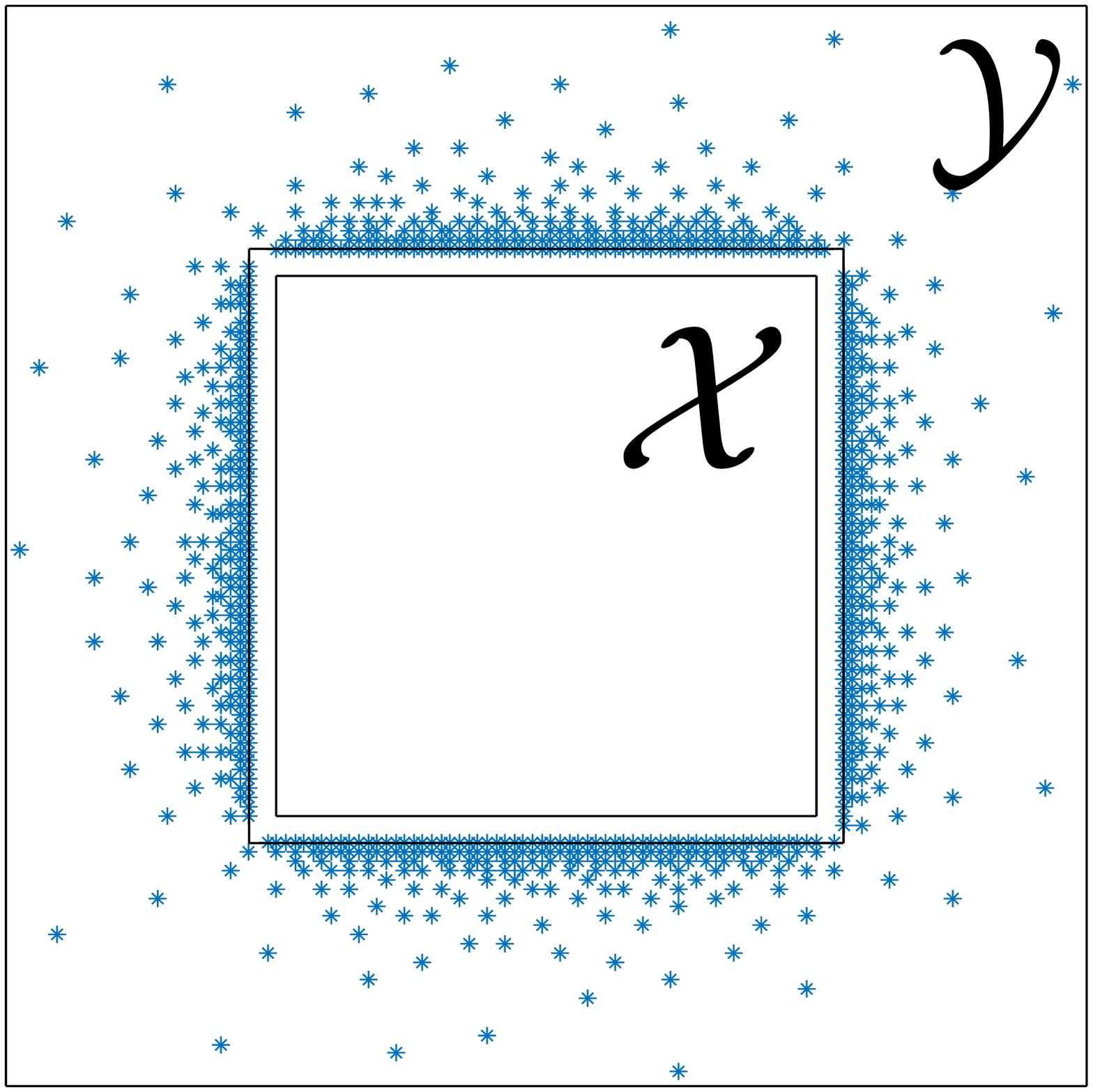}}
	\hspace{1em}
	\subfloat[$\frac{1}{|X_0|-|X_\text{rep}|}\sum_{x_i\in X_0\backslash X_\text{rep}}|e_i(y)|$ for 
	${y \in \mathcal{Y}}$ in $\log_{10}$ scale.]
	{\includegraphics[height=0.23\textheight]{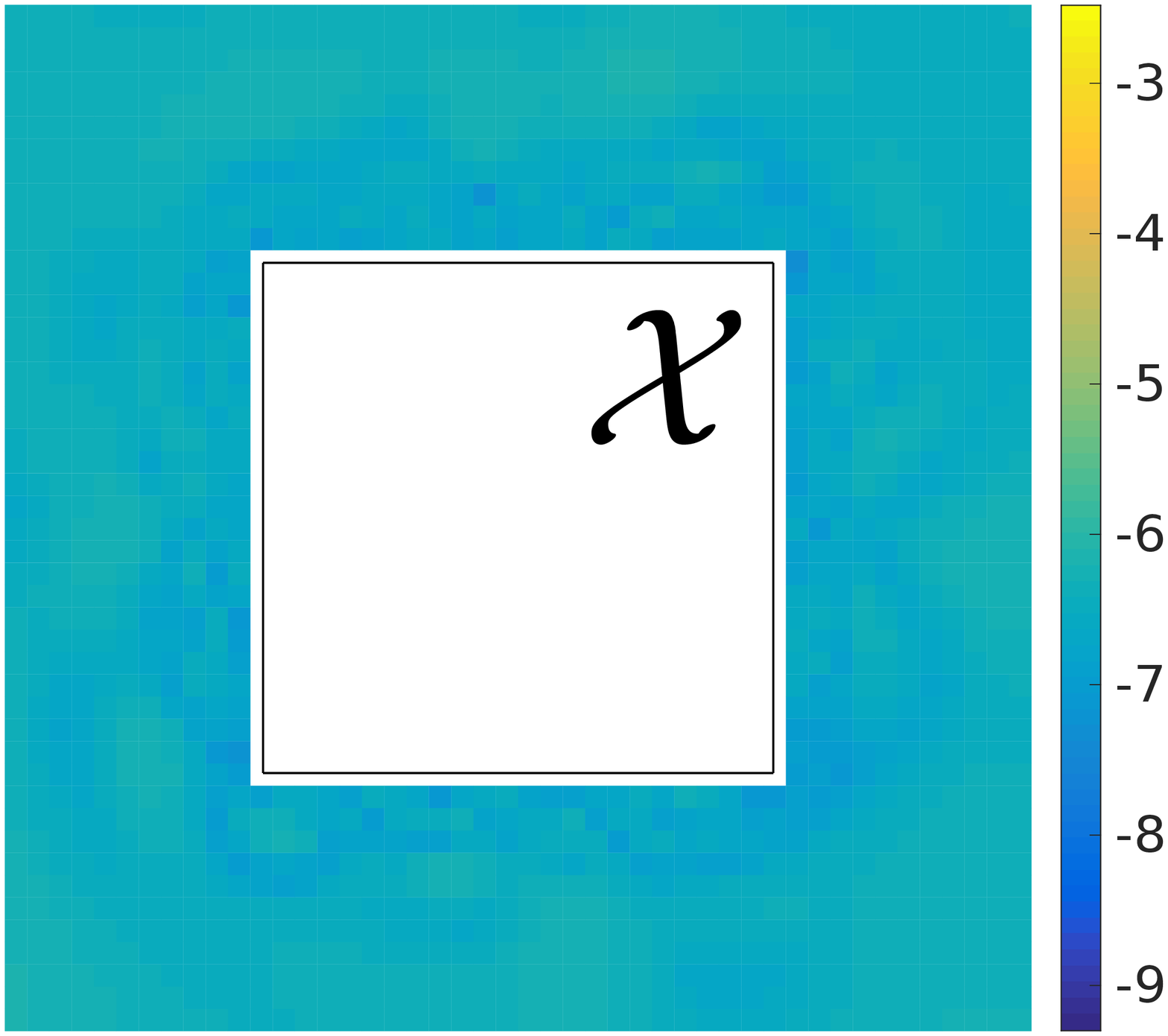}}
	\caption{$Y_p^\text{id}$ obtained by \cref{alg:proxy} with adaptively selected initial grid
		$Y_d$ and the associated error distribution in $\mathcal{Y}$, for 
		$K_1(x,y) = \frac{1}{|x-y|}$ with the nearby domain case in 2D.
		The plotted area is $[-2,2]^2$. The error function values 
		outside of $[-2,2]^2$ are approximately $10^{-7}$.
		The maximum average error over $\mathcal{Y}$ is found to be $7.97\times 10^{-7}$.
		\label{fig:test4_weak_adaptive}}
\end{figure}

A heuristic solution is to adaptively select more candidate points in the area where $K(x,y)$ has larger
variation (e.g., according to $|\nabla_y K(x,y)|$). 
To test this idea, we uniformly select half of the candidate points, approximately 7500 points,
from the small area $[2,2]^2\backslash[1.1,1.1]^2$  and the other half from the remaining large area $[9,9]^2\backslash[2,2]^2$. 
The corresponding results are shown in \cref{fig:test4_weak_adaptive}. 
More proxy points ($|Y_p^\text{id}| = 909$) are selected, especially from the area near $\mathcal{X}$, 
and the obtained average error also meets the expected accuracy at any $y\in \mathcal{Y}$.  
These results corroborate our explanation of the possible cause of the accuracy degradation. 
Alternatively, we can consider re-applying \cref{alg:proxy} with denser initial candidates in areas with 
large error to improve the quality of proxy point set $Y_p^\text{id}$. 

Note that in the weak admissibility setting, 
%if $\mathcal{Y}$ is defined as $[-9,9]^d\backslash [-1,1]^d$ to meet the weak admissibility condition with $\mathcal{X}$, 
$K_1(x,y)$ is singular on the boundary between $\mathcal{X}$ and $\mathcal{Y}$ and no KL expansion exists for 
$K_1(x,y)$ over $\mathcal{X}\times \mathcal{Y}$. In this case, \cref{alg:proxy} and the proposed 
algorithm no longer work. Practically, we can add a small gap between $\mathcal{X}$ and $\mathcal{Y}$ but 
$r_\text{KL}$ would be large and numerical tests show that large numbers of points for $Y_d$ and $Y_p^\text{id}$
are needed to achieve the same accuracy. 

\begin{figure}
	\centering
	\includegraphics[width=\textwidth]{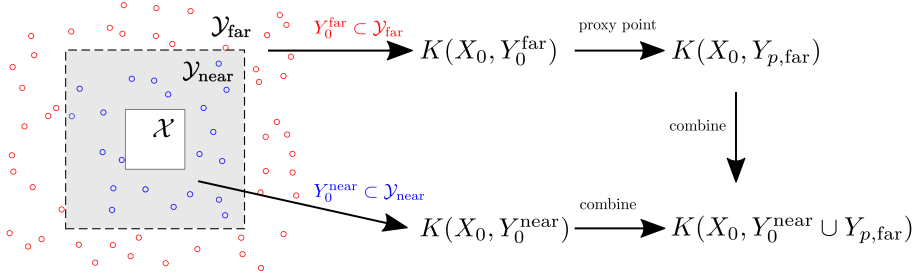}
	\caption{Hybrid variant of the proposed algorithm. Approximation of target matrix $K(X_0, Y_0)$ is
		obtained by an ID approximation of a smaller matrix $K(X_0, Y_0^\text{near}\cup Y_{p,\text{far}})$
		where $|Y_{p,\text{far}}|$ is not related to $|Y_0|$ and $|Y_0^\text{near}|$ is usually expected to 
		be nearly constant in real problems. \label{fig:hybrid}		
	}
\end{figure}

Another solution to both the accuracy degradation and the kernel singularity issue is inspired by the hybrid
method in \cite{ho_fast_2012} and is illustrated in \cref{fig:hybrid}. 
Consider $K_1(x,y)$ with the domain pair $\mathcal{X}\times \mathcal{Y}$ in 2D that satisfies the weak 
admissibility condition.
Split $\mathcal{Y}$ into a neighboring field $\mathcal{Y}_\text{near}$ and a far field $\mathcal{Y}_\text{far}$.
From the previous tests, the proposed algorithm works well for 
$\mathcal{X}\times \mathcal{Y}_\text{far}$ but does not work or has large error for $\mathcal{X}\times\mathcal{Y}_\text{near}$. 
For a target point set $Y_0\subset \mathcal{Y}$, split it as $Y_0 = Y_0^\text{near} \cup Y_0^\text{far}$ so that 
$Y_0^\text{near} \subset \mathcal{Y}_\text{near}$ and $Y_0^\text{far} \subset \mathcal{Y}_\text{far}$. 
The idea is to only apply the proposed algorithm over $K_1(X_0, Y_0^\text{far})$ and directly work on 
$K_1(X_0, Y_0^\text{near})$. Specifically, denote $Y_{p,\text{far}}$ as some proxy point set selected for 
$K_1(x,y)$ over $\mathcal{X}\times \mathcal{Y}_\text{far}$ and find an ID approximation of 
$K_1(X_0, Y_0^\text{near}\cup Y_{p,\text{far}})$ by sRRQR as
\[
K_1(X_0, Y_0^\text{near}\cup Y_{p,\text{far}}) \approx W_\text{rep} 
K_1(X_\text{rep}, Y_0^\text{near}\cup Y_{p,\text{far}}).
\]
The ID approximation of $K_1(X_0,Y_0)$ is then defined as $W_\text{rep} K_1(X_\text{rep}, Y_0)$. 
In general, the splitting of $\mathcal{Y}$ should be kernel-dependent. 
This hybrid method will be illustrated in the next subsection on $\mathcal{H}^2$ construction. 

\subsection{$\mathcal{H}^2$ matrix construction}
We now consider the ID-based $\mathcal{H}^2$ construction of symmetric kernel matrices $K(X, X)$ with some 
prescribed point set $X \subset \mathbb{R}^d$. The following two admissibility conditions are considered. 
\begin{itemize}
\item strong admissibility condition: For any two non-adjacent boxes, the two enclosed point clusters are admissible.
\item weak admissibility condition: For any two non-overlapping boxes, the two enclosed point clusters are admissible. 
The $\mathcal{H}^2$ matrix with this admissibility condition is usually called an HSS matrix in the literature. 
\end{itemize}
The associated $\mathcal{H}^2$ representations are referred to as 
$\mathcal{H}^2_\text{strong}$ and $\mathcal{H}^2_\text{weak}$ respectively.

For problems in $\mathbb{R}^d$, to maintain constant point density, $N$ points are uniformly and randomly 
distributed in a box with equal edge length $L = N^{1/d}$. 
At the $k$th level of the hierarchical partitioning of the box, the sub-domains contain $2^{dk}$ boxes with edge
length $L/2^{k}$. 
For the associated \cref{problem:targetmatrix} at the $k$th level, set $\mathcal{X} = [-\frac{L}{2^{k+1}}, \frac{L}{2^{k+1}}]^d$ and
\begin{equation*}
\mathcal{Y} = 
\left\{\begin{array}{ll}
[-(L-\frac{L}{2^{k+1}}), L-\frac{L}{2^{k+1}}]^d\setminus [-\frac{3L}{2^{k+1}}, \frac{3L}{2^{k+1}}]^d & \text{  for } \mathcal{H}^2_\text{strong} \\ 

[-(L-\frac{L}{2^{k+1}}), L-\frac{L}{2^{k+1}}]^d\setminus \mathcal{X} & \text{  for } \mathcal{H}^2_\text{weak}
\end{array}\right.,
\end{equation*}
based on the admissibility conditions. 

For $\mathcal{H}^2_\text{strong}$ construction, the proposed algorithm with both \textit{ID} and \textit{Random Selection}
schemes for $Y_p$ is compared with the ID approximation using sRRQR. $Y_p^\text{id}$ is selected by \cref{alg:proxy} from 
a uniform initial set pair $X_d\times Y_d$ with approximately $1500\times 15000$ points and $Y_p^\text{rand}$ contains
2000 points in $\mathcal{Y}$. 

For $\mathcal{H}^2_\text{weak}$ construction, the hybrid algorithm in \cref{fig:hybrid} with both {\it Random} and 
{\it ID Selection} schemes for $Y_{p,\text{far}}$ (denoted as $Y_{p,\text{far}}^\text{rand}$ and $Y_{p,\text{far}}^\text{id}$)
is also tested. $\mathcal{Y}_\text{far}$, $Y_{p,\text{far}}^\text{rand}$ and $Y_{p,\text{far}}^\text{id}$ used in 
this algorithm are the same as $\mathcal{Y}$, $Y_p^\text{rand}$ and $Y_p^\text{id}$ in the $\mathcal{H}^2_\text{strong}$ 
construction case above.
The non-hybrid version of the proposed algorithm for $K_1(x,y)$ is not tested due to the singularity of $K_1(x,y)$ on the boundary 
between $\mathcal{X}$ and $\mathcal{Y}$ in the weak admissibility setting. 
%With the defined $\mathcal{X}\times\mathcal{Y}$  for $\mathcal{H}^2_\text{weak}$, a degradation of accuracy when $Y_p^\text{id}$
%is selected from the uniform candidate set $Y_d$ has also been observed for $K_2(x,y)$ in some simple tests not shown here. 
To select $Y_p^\text{id}$ for $K_2(x,y)$ by \cref{alg:proxy}, following the strategy from the previous subsection, half of the candidate 
points $Y_d$ are uniformly selected from a small area near $\mathcal{X}$ with the other half from the remaining large area. 
A similar strategy was used for $Y_p^\text{rand}$.  

\begin{figure}[h]
	\centering
	\subfloat[$K_1(x,y)=1/|x-y|$]
	{\includegraphics[width=0.45\textwidth]{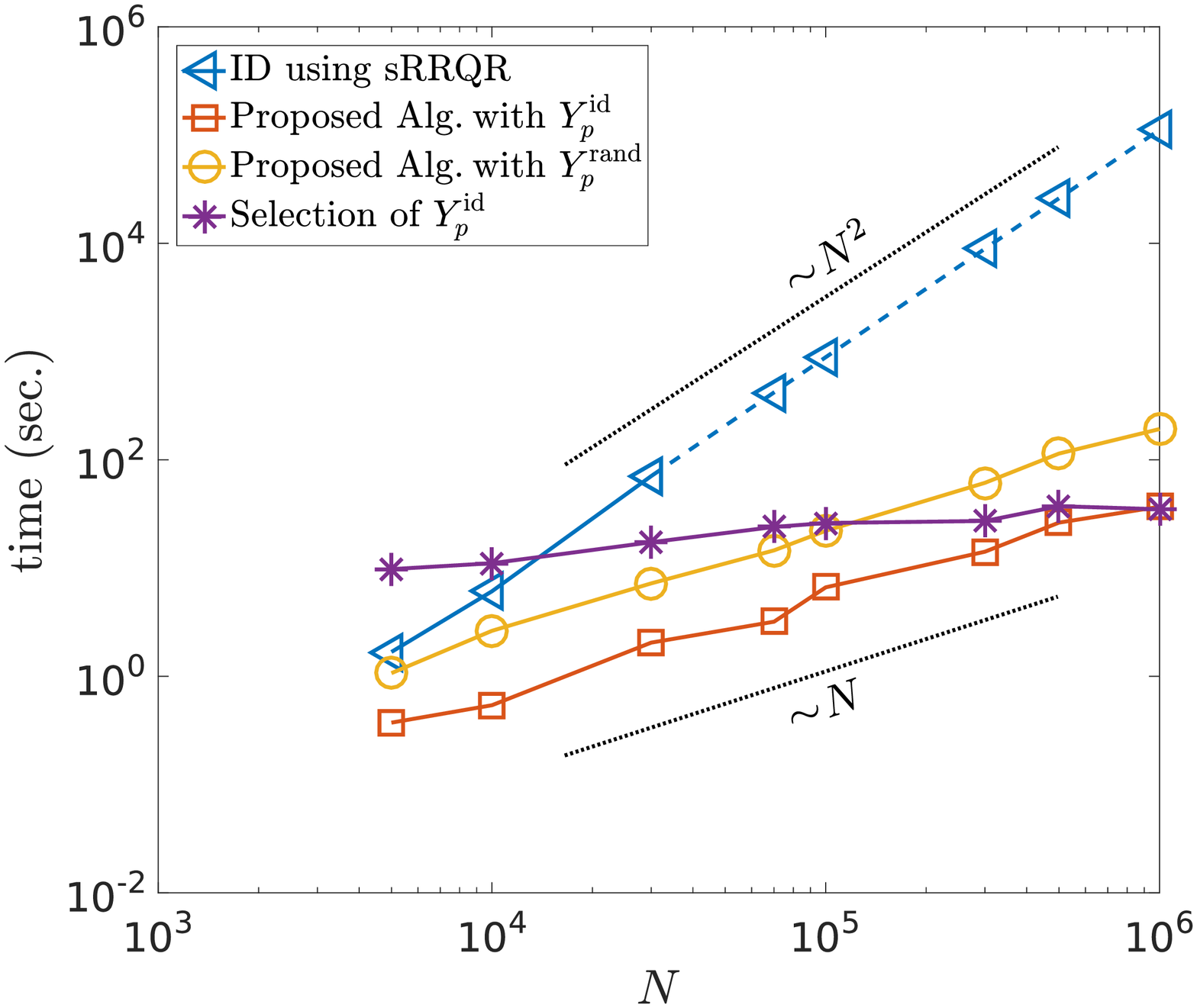}}
	\hspace{1em}
	\subfloat[$K_2(x,y)=\sqrt{1+|x-y|^2}$]
	{\includegraphics[width=0.45\textwidth]{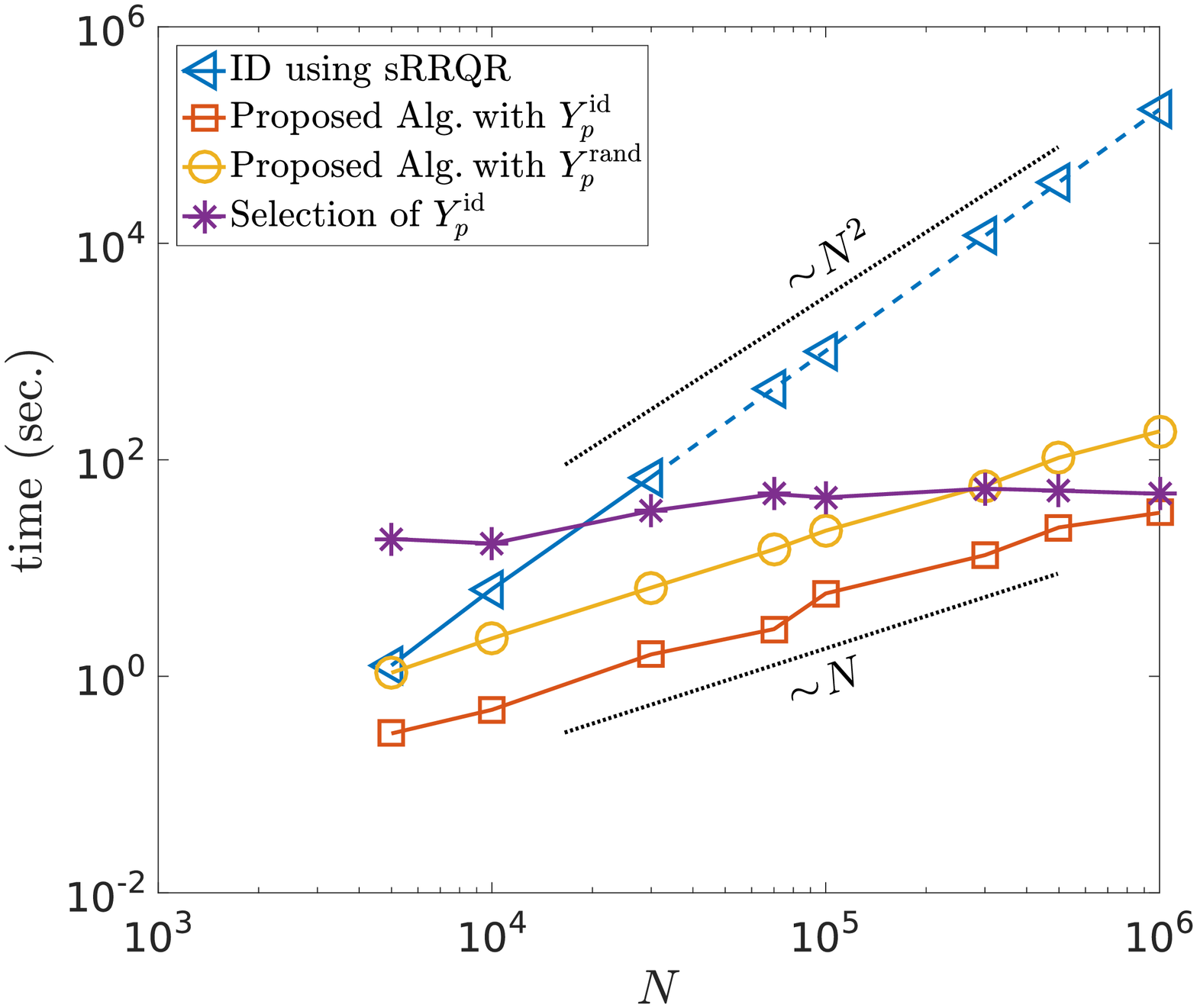}}
	\caption{Runtime of $\mathcal{H}^2_\text{strong}$ construction for two kernels in 2D. 
	Construction with ID using sRRQR is not tested when $N \geqslant 70000$ due to 
	memory limitations and the dashed lines indicate an extrapolation of the data. 
	\label{fig:test5_H2_time}}
\end{figure}

\begin{figure}[h]
	\centering
	\subfloat[$K_1(x,y)=1/|x-y|$]
	{\includegraphics[width=0.45\textwidth]{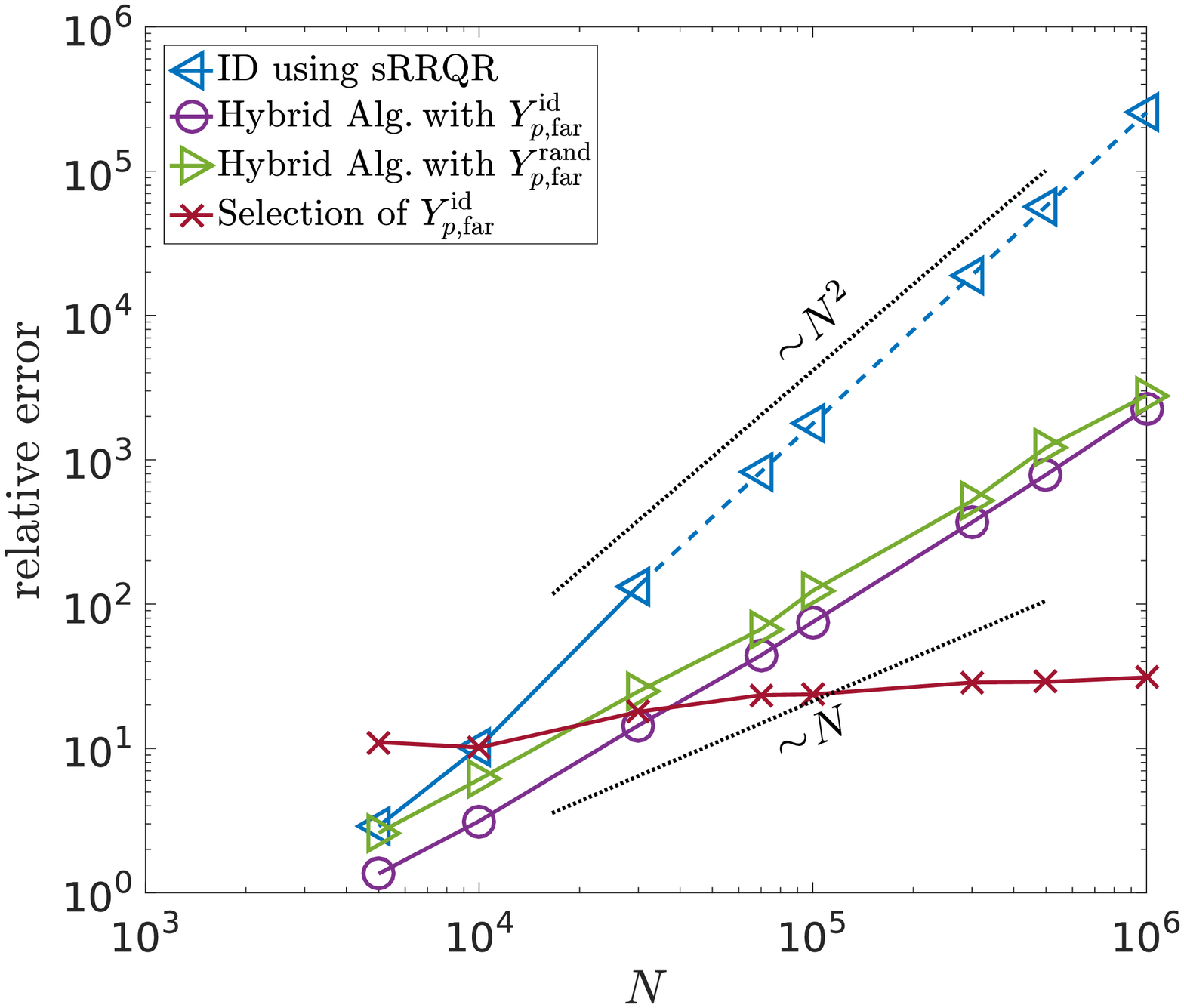}}
	\hspace{1em}
	\subfloat[$K_2(x,y)=\sqrt{1+|x-y|^2}$]
	{\includegraphics[width=0.45\textwidth]{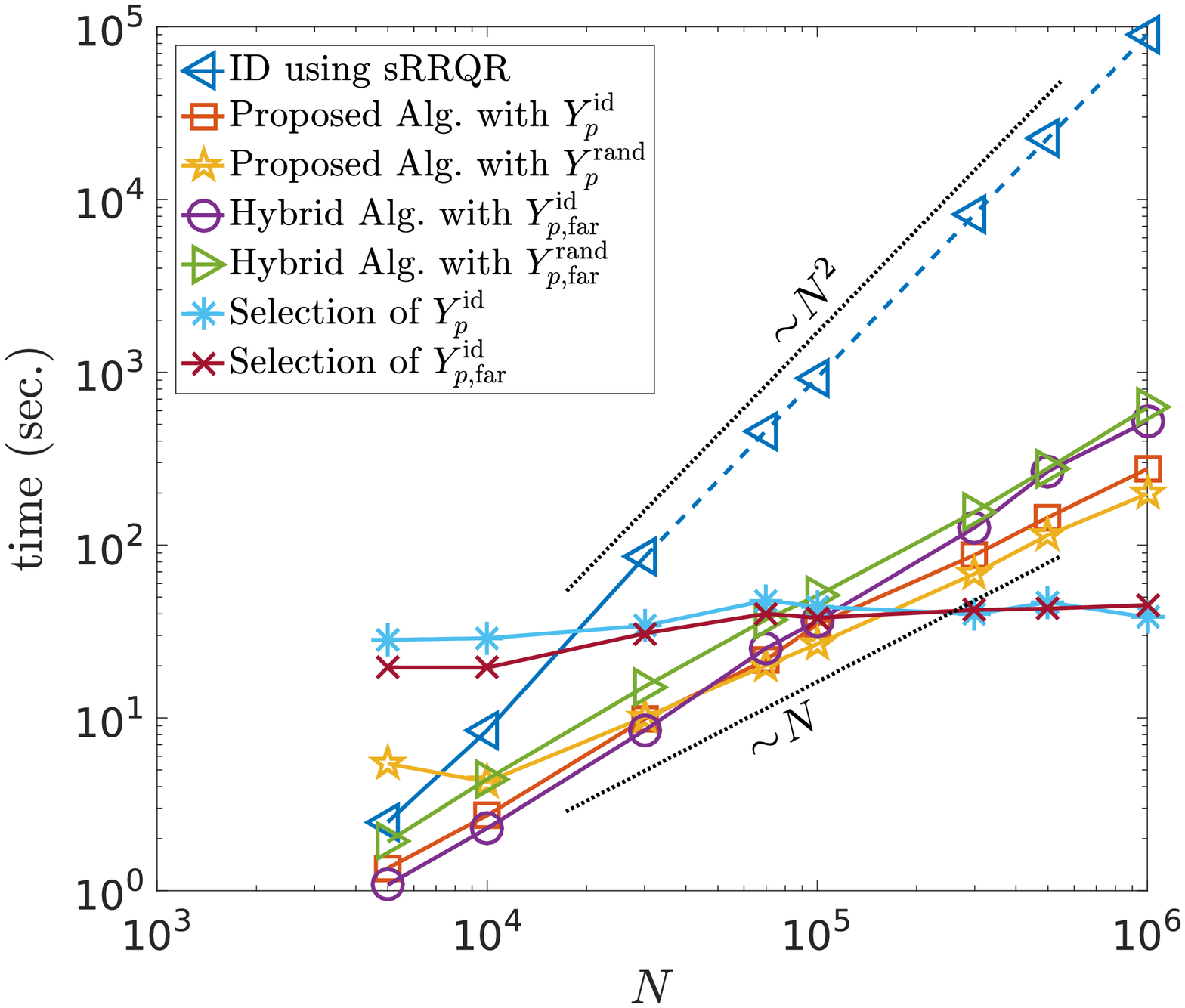}}
	\caption{Runtime of $\mathcal{H}^2_\text{weak}$ construction for two kernels in 2D.
		The dashed lines indicate an extrapolation of the data. 
	\label{fig:test5_HSS_time}}
\end{figure}

In the hierarchical partitioning, a sub-domain is subdivided when it has more than 300 points. 
For all the ID approximations in $\mathcal{H}^2$
construction, a relative error threshold of $\tau = 10^{-6}$ is applied. We consider the two prescribed kernels in 2D.
\cref{fig:test5_H2_time} and \cref{fig:test5_HSS_time} show the runtime of our sequential
Matlab implementation for $\mathcal{H}^2_\text{strong}$ and $\mathcal{H}^2_\text{weak}$ construction.
\cref{tab:test5_H2} and \cref{tab:test5_HSS} list some detailed data of these constructions.

\begin{table}[ht]
\centering
\caption{Numerical results for $\mathcal{H}^2_\text{strong}$ construction. 
$\mathcal{S}$ denotes the storage cost (MB) of the obtained matrix representation. $\mathcal{E}$ denotes
the relative error of the approximations of all compressed blocks. $\mathcal{S}$ or $\mathcal{E}$ with 
subscripts \textnormal{sRRQR}, $Y_p^\text{id}$ and $Y_p^\text{rand}$ are the associated values 
obtained by the $\mathcal{H}^2$ construction using \textnormal{ID approximation using sRRQR, the 
proposed algorithm with} $Y_p^\text{id}$
and that with $Y_p^\text{rand}$, respectively. $\mathcal{S}_\text{N}$ denotes the storage cost of 
the original dense matrix and $\mathcal{S}_\text{inadm}$ denotes the storage cost of the dense
blocks in the $\mathcal{H}^2$ representation.\label{tab:test5_H2} }
\begin{tabular}{c|ccccc|ccc}
\toprule
		  \multicolumn{8}{c}{$K_1(x,y) = 1 / |x-y|$}   \\
$N$	 & $\mathcal{S}_N$ & $\mathcal{S}_\text{inadm}$ & $\mathcal{S}_\text{sRRQR}$ & 
$\mathcal{S}_{Y_p^\text{id}}$ & $\mathcal{S}_{Y_p^\text{rand}}$ & $\mathcal{E}_\text{sRRQR}$ & 
$\mathcal{E}_{Y_p^\text{id}}$ & $\mathcal{E}_{Y_p^\text{rand}}$ \\
\midrule
5e3  & 1.9e2  & 25     & 35    & 39    &  37    & 6.9e-7 & 1.1e-6 & 1.5e-6\\
1e4  & 7.6e2  & 90     & 1.0e2 & 1.1e2 &  1.1e2 & 7.9e-7 & 1.2e-6 & 2.0e-6\\
3e4  & 6.9e3  & 2.2e2  & 2.9e2 & 3.1e2 &  2.7e2 & 9.1e-7 & 1.4e-6 & 5.7e-6\\
7e4  & 3.7e4  & 1.1e3  & -     & 1.3e3 &  1.2e3 & -      & 1.7e-6 & 1.9e-5\\
1e5  & 7.6e4  & 6.4e2  & -     & 9.9e2 &  8.1e2 & -      & 1.8e-6 & 6.9e-5\\
3e5  & 6.9e5  & 3.8e3  & -     & 4.5e3 &  4.2e3 & -      & 4.1e-6 & 1.4e-4\\
5e5  & 1.9e6  & 4.1e3  & -     & 5.3e3 &  4.6e3 & -      & 8.0e-6 & 2.1e-4\\
1e6  & 7.6e6  & 1.6e4  & -     & 1.8e4 &  1.7e4 & -      & 9.5e-6 & 2.0e-4\\

\bottomrule
\end{tabular}

\vspace{0.7em}
\begin{tabular}{c|ccccc|ccc}
\toprule
		  \multicolumn{8}{c}{$K_2(x,y) = \sqrt{1+|x-y|^2}$}   \\
$N$	 & $\mathcal{S}_N$ & $\mathcal{S}_\text{inadm}$ & $\mathcal{S}_\text{sRRQR}$ & 
$\mathcal{S}_{Y_p^\text{id}}$ & $\mathcal{S}_{Y_p^\text{rand}}$ & $\mathcal{E}_\text{sRRQR}$ & 
$\mathcal{E}_{Y_p^\text{id}}$ & $\mathcal{E}_{Y_p^\text{rand}}$ \\
\midrule
5e3  & 1.9e2  & 25     & 31    & 33    &  30    & 8.1e-7 & 4.4e-7 & 1.5e-6\\
1e4  & 7.6e2  & 90     & 97    & 1.0e2 &  97    & 8.8e-7 & 4.4e-7 & 1.2e-6\\
3e4  & 6.9e3  & 2.2e2  & 2.5e2 & 2.7e2 &  2.5e2 & 1.1e-6 & 4.7e-7 & 2.9e-6\\
7e4  & 3.7e4  & 1.1e3  & -     & 1.2e3 &  1.2e3 & -      & 5.0e-7 & 3.2e-6\\
1e5  & 7.6e4  & 6.4e2  & -     & 8.0e2 &  7.1e2 & -      & 4.9e-7 & 3.5e-6\\
3e5  & 6.9e5  & 3.8e3  & -     & 4.2e3 &  3.9e3 & -      & 4.7e-7 & 4.6e-6\\
5e5  & 1.9e6  & 4.1e3  & -     & 4.7e3 &  4.3e3 & -      & 5.1e-7 & 5.3e-6\\
1e6  & 7.6e6  & 1.6e4  & -     & 1.7e4 &  1.7e4 & -      & 5.1e-7 & 5.5e-6\\
\bottomrule
\end{tabular}
\end{table}

\begin{table}[ht]
\centering
\caption{Numerical results for $\mathcal{H}^2_\text{weak}$ construction. Subscripts 
$Y_{p,\text{far}}^\text{id}$ and $Y_{p,\text{far}}^\text{rand}$ correspond to the hybrid
algorithm with $Y_{p,\text{far}}^\text{id}$ and $Y_{p,\text{far}}^\text{rand}$ as the 
proxy point set for $\mathcal{Y}_\text{far}$ respectively.
See the caption for \cref{tab:test5_H2} for the definitions of other notations.
\label{tab:test5_HSS} }
\begin{tabular}{c|ccccc|ccc}
\toprule
		  \multicolumn{8}{c}{$K_1(x,y) = 1 / |x-y|$}   \\
$N$	 & $\mathcal{S}_N$ & $\mathcal{S}_\text{inadm}$ & $\mathcal{S}_\text{sRRQR}$ & 
$\mathcal{S}_{Y_{p,\text{far}}^\text{id}}$ & $\mathcal{S}_{Y_{p,\text{far}}^\text{rand}}$ & $\mathcal{E}_\text{sRRQR}$ & 
$\mathcal{E}_{Y_{p,\text{far}}^\text{id}}$ & $\mathcal{E}_{Y_{p,\text{far}}^\text{rand}}$ \\
\midrule
5e3 & 1.9e2  & 3.9    & 28    & 29    &  29    & 6.2e-6 & 1.8e-6 & 5.7e-7\\
1e4 & 7.6e2  & 12     & 67    & 69    &  69    & 8.9e-6 & 3.4e-6 & 1.1e-6\\
3e4 & 6.9e3  & 27     & 2.7e2 & 2.7e2 &  2.7e2 & 1.1e-5 & 3.8e-6 & 1.7e-6\\
7e4 & 3.7e4  & 1.4e2  & -     & 7.9e2 &  7.9e2 & -      & 4.6e-6 & 2.6e-6\\
1e5 & 7.6e4  & 75     & -     & 1.1e3 &  1.2e3 & -      & 5.5e-6 & 3.0e-6\\
3e5 & 6.9e5  & 4.8e2  & -     & 4.3e3 &  4.3e3 & -      & 7.1e-6 & 4.6e-6\\
5e5 & 1.9e6  & 4.7e2  & -     & 7.6e3 &  7.6e3 & -      & 1.0e-5 & 7.2e-6\\
1e6 & 7.6e6  & 1.9e3  & -     & 1.7e4 &  1.7e4 & -      & 1.3e-5 & 9.4e-6\\

\bottomrule
\end{tabular}

\vspace{0.7em}
{\setlength \tabcolsep{3pt}
\begin{tabular}{c|ccccc|ccccc}
\toprule
		  \multicolumn{11}{c}{$K_2(x,y) = \sqrt{1+|x-y|^2} $}   \\
$N$	 &  $^\ast\mathcal{S}_\text{sRRQR}$ & 
$\mathcal{S}_{Y_p^\text{id}}$ & $\mathcal{S}_{Y_p^\text{rand}}$ & 
$\mathcal{S}_{Y_{p,\text{far}}^\text{id}}$ & $\mathcal{S}_{Y_{p,\text{far}}^\text{rand}}$ &
$\mathcal{E}_\text{sRRQR}$ & 
$\mathcal{E}_{Y_p^\text{id}}$ & $\mathcal{E}_{Y_p^\text{rand}}$ &
$\mathcal{E}_{Y_{p,\text{far}}^\text{id}}$ & $\mathcal{E}_{Y_{p,\text{far}}^\text{rand}}$ \\
\midrule
5e3  & 15    & 37    &  23    & 17    &  15    & 1.7e-6 & 6.5e-7 & 1.4e-6 & 1.3e-6 & 3.3e-6\\
1e4  & 33    & 86    &  47    & 41    &  34    & 1.9e-6 & 8.2e-7 & 3.5e-6 & 1.2e-6 & 3.0e-6\\
3e4  & 94    & 3.0e2 &  1.0e2 & 1.4e2 &  1.0e2 & 2.6e-6 & 8.7e-7 & 7.0e-6 & 1.7e-6 & 3.4e-6\\
7e4  & -     & 7.6e2 &  2.6e2 & 4.0e2 &  3.3e2 & -      & 1.0e-6 & 1.3e-5 & 2.1e-6 & 2.9e-6\\
1e5  & -     & 1.0e3 &  2.2e2 & 5.0e2 &  3.5e2 & -      & 1.1e-6 & 1.5e-5 & 2.1e-6 & 3.5e-6\\
3e5  & -     & 3.1e3 &  7.6e2 & 1.8e3 &  1.3e3 & -      & 4.5e-6 & 2.1e-5 & 2.7e-6 & 3.4e-6\\
5e5  & -     & 4.9e3 &  8.0e2 & 2.7e3 &  1.7e3 & -      & 3.6e-6 & 2.3e-5 & 2.7e-6 & 3.7e-6\\
1e6  & -     & 9.9e3 &  2.4e3 & 6.3e3 &  4.3e3 & -      & 5.4e-6 & 2.4e-5 & 3.3e-6 & 3.5e-6\\
\bottomrule
\end{tabular}}
{\footnotesize $^*$Refer to the table for $K_1(x,y)$ for values of $\mathcal{S}_N$ and $\mathcal{S}_\text{inadm}$.\hfill }
\end{table}

For both constructions, the runtime of the pre-calculation for $Y_p^\text{id}$ is significant but 
its asymptotic complexity is only $O(\log(N))$ as the selection of $Y_p^\text{id}$ is only performed
once for each level. 
Both the proposed and hybrid algorithms lead to nearly linear $\mathcal{H}^2$ 
construction which can also be verified by complexity analysis if we assume that the maximum rank of 
the low-rank off-diagonal blocks and the number of points in $Y_p$ for each level are both of constant scale. 
Also, $\mathcal{H}^2$ construction with these two algorithms has larger storage cost compared to 
that with ID using sRRQR because, as explained earlier, they generally select more rows in the ID
approximation.

For $\mathcal{H}^2_\text{strong}$ construction, the proposed algorithm with $Y_p^\text{rand}$ takes more
time since $Y_p^\text{rand}$ has more points than $Y_p^\text{id}$  which contains approximately
900 points for each level. 
%Consistent with the previous tests, $Y_p^\text{id}$ here gives a better approximation compared to 
%$Y_p^\text{rand}$.
The relative errors of $\mathcal{H}^2_\text{strong}$ and $\mathcal{H}^2_\text{weak}$ approximations for 
$K_1(x,y)$ both increase with larger $N$ which can be observed for both proxy point selection schemes and 
for both the proposed algorithm and the ID using sRRQR. 
This is mainly due to the amplification of errors at the level-by-level $\mathcal{H}^2$ construction
and is also kernel-dependent. The hierarchical partitioning trees have $4, 4, 5, 6, 6, 7, 7, 7$ levels
for the values of $N$ tested, which roughly matches the incremental pattern of the errors in \cref{tab:test5_H2}. 

For $\mathcal{H}^2_\text{weak}$ construction, the hybrid algorithms with both selection schemes for $Y_p$
are effective and provide good approximations. In addition, for $K_2(x,y)$, the hybrid algorithm
has less storage cost (i.e., smaller ranks for ID approximations) and similar or even smaller relative
errors compared to the non-hybrid algorithm with $Y_p^\text{id}$. This advantage of the hybrid algorithm 
is expected since the hybrid algorithm directly works on a part of $Y_0$ in the ID approximation.

\section{Conclusion}
We proposed an efficient low-rank approximation algorithm for the sub-blocks of kernel matrices 
that can also be regarded as a generalization of the proxy-surface method.  
For the proposed algorithm, two proxy point selection schemes were introduced in \cref{sec:proxy_point}
as well as two heuristic improvements to the schemes in \cref{sec:improvement}. 
The two proxy point selection schemes that were introduced are general in that they can be applied to any 
kernel and any domain pair. It should be possible to design specialized selection schemes that are kernel-dependent
and thus are more effective than the general schemes, as long as condition \cref{prop:Ysample} is met. 
%It is worth noting that the selection of the proxy points is flexible in the proposed algorithm, 
%kernel-dependent and more effective proxy point selection schemes may exist. 
%For a block $K(X_0, Y_0)$ with $X_0\times Y_0$ at a specific domain pair $\mathcal{X}\times\mathcal{Y}$, 
%the algorithm has complexity $O(r|X_0|)$, independent of $|Y_0|$.
In practice, the algorithm can be used for hierarchical matrix construction for general translation-invariant
kernels to give a construction cost linear in the matrix dimension if the maximum rank of the low-rank off-diagonal 
blocks does not increase with the matrix dimension.

\bibliographystyle{plain}
\bibliography{references}

\end{document}